\documentclass[12pt]{article}
\usepackage{amssymb,a4}
\usepackage{mathrsfs}
\usepackage{epsf}
\usepackage{bm}
\usepackage{latexsym}
\usepackage{amsmath,amsfonts,amssymb,amsthm}
\RequirePackage{amsopn} \RequirePackage{amsfonts}
\RequirePackage{amsthm}

\newcommand{\1}{{{\bf 1}}}

\newcommand{\Z}{\mathbb{Z}}

\newcommand{\C}{\mathbb{C}}
\newcommand{\N}{\mathbb{N}}

\def\C{{\mathbb C}}

\def\R{{\mathbb R}}

\def\Z{{\mathbb Z}}

\def\N{{\mathbb N}}
\def\1{{\bf 1}}

\def \pf{\noindent {\bf Proof: \,}}
\def\theequation{5.\arabic{equation}}
\def \a{\alpha}
\def \b{\beta}
\def \e{\epsilon}
\def \h{\mathfrak{h}}
\def \g{\mathfrak{g}}
\def \l{\lambda}
\def \w{\omega}
\begin{document}
\newtheorem{theorem}{Theorem}[section]
\newtheorem{proposition}[theorem]{Proposition}
\newtheorem{lem}[theorem]{Lemma}
\newtheorem{corollary}[theorem]{Corollary}
\newtheorem{definition}[theorem]{Definition}
\newtheorem{example}[theorem]{Example}
\newtheorem{remark}[theorem]{Remark}

\numberwithin{equation}{section}

\newenvironment{namelist}[1]{%
\begin{list}{}
{
\let\makelabel\namelistlabel
\settowidth{\labelwidth}{#1}
 \setlength{\leftmargin}{1.1\labelwidth}
 }
 }{%
 \end{list}}
 \begin{center}
{\Large {\bf Unitary vertex operator algebras}} \\
\vspace{0.5cm} \vspace{0.5cm} Chongying Dong
\\
Department of Mathematics, Sichuan University\\
 Chengdu China \&\\
Department of Mathematics, University of
California\\ Santa Cruz, CA 95064 \\

\vspace{0.5cm}
Xingjun Lin\\
Department of Mathematics, Sichuan University\\
 Chengdu China\\

\end{center}

\begin{abstract}
Unitary vertex operator algebras are introduced and studied. It is proved
that most well-known rational vertex operator algebras are unitary. The classification of unitary vertex operator algebras with central charge $c\leq 1$ is also discussed.
\end{abstract}
\section{Introduction \label{intro}}
\def\theequation{1.\arabic{equation}}
\setcounter{equation}{0}

Both bilinear form and Hermitian form are power tools in the study of general algebras and their representations.
Based on a symmetric contravariant bilinear form with respect to a Cartan involution in \cite{B}  for vertex algebras associated to even lattices, notions of invariant bilinear form for an vertex operator algebra and contragradient module were introduced and studied in \cite{FHL}. The space of invariant bilinear form for an arbitrary vertex operator algebra was determined in \cite{Li}. In this paper, we study the invariant Hermitian form on vertex
operator algebras and their modules, and use the Hermitian form to define and investigate
the unitary vertex operator algebras.

One important motivation for studying the unitary vertex operator algebra comes from the unitary representations of infinite dimensional Lie algebras such as Virasoro algebras and Kac-Moody algebras. The unitary highest weight representations of these
algebras produce fundamental families of rational vertex operator algebras. It is well known today that the theory of vertex operator algebra unifies representation theory of many infinite dimensional Lie algebra via locality.
So it is natural to have a notion  of unitary vertex operator algebra so that in the case of Virasoro and affine Kac-Moody algebras, these two unitarities are equivalent. While the  Hermitian form in the the theory of infinite dimensional Lie
algebras is defined to be the Hermitian form which is contravariant under some anti-linear anti-involution of
corresponding universal enveloping algebras, the Hermitian forms on vertex operator algebras and their modules are defined to be contravariant under some anti-linear involutions. In the case of Virasoro and affine vertex operator algebras, the anti-linear involutions induce anti-linear anti-involutions which are exactly the given ones of the infinite dimensional Lie algebras
which are required for the unitary representations. The unitary vertex operator algebra is defined to be a vertex operator algebra associated with a positive definite Hermitian form.  We will prove that the vertex operator algebras associated
to the unitary highest weight representations for the Heisenberg algebra, Virasoro algebra and affine Kac-Moody algebras are indeed the unitary vertex operator algebras. We also discuss the unitarity of irreducible modules for these vertex operator algebras.

The positive definite Hermitian form already appeared in
the theories of lattice vertex operator algebras \cite{B},
\cite{FLM}, and their $\Z_2$-orbifold vertex operator algebras \cite{DGM}.
We show that the lattice vertex operator algebras
associated to positive definite even lattices are unitary vertex
operator algebras. We also establish the unitarity for the irreducible modules and the irreducible $\theta$-twisted modules for lattice vertex operator algebras where $\theta$ is the automorphism of lattice vertex operator algebras induced from the $-1$ isometry of the lattices. These results are then used to
show that the moonshine vertex operator algebra $V^{\natural}$ is also unitary.

Another motivation comes from connection between algebraic and analytic approaches to 2-dimensional conformal field theory.
While the  algebraic approach uses vertex operator algebras, the analytic approach uses conformal nets. It has been expected
that these two approaches are equivalent in the following sense: one can construct conformal nets and vertex operator algebras from each other. Although it is not clear how this can be achieved, one can see the similarity of these two approaches in many examples. The basic object in the theory of operator algebras is Hilbert space. So it is desearble to have a positive definite Hermitian form on vertex operator algebra whose completion give rise to a Hilbert space. From this point of view, studying the unitary vertex operator algebra is  the first step in constructing conformal nets from vertex operator algebras.

This paper is organized as follows. In Section 2, we introduce
the notion of unitary vertex operator algebra, and give some
elementary facts about unitary vertex operator algebras. In
Section 3, we show that the unitary rational and $C_2$-cofinite vertex
operator algebras could be extended to a unitary vertex operator algebra by a simple current under some
assumption. In Section 4, we prove that some well-known vertex
operator algebras are unitary. In Section 5, we give some
results about the classification of unitary vertex operator
algebras with central charge $c\leq 1$.
\section{Preliminaries }\def\theequation{2.\arabic{equation}}
\setcounter{equation}{0}
We assume that the readers are familiar with the notion of vertex operator algebra and the basic facts about vertex operator algebra as presented in \cite{FLM}, \cite{FHL}, \cite{DLM1}, \cite{DLM2}, \cite{LL} and \cite{Z}. In this paper, we only consider the vertex operator algebra $(V, Y, \1, \w)$ of CFT-type, i.e. $V_n=0, n<0$ and $V_0=\mathbb{C}\bf 1$.
\begin{definition}{\rm
Let $(V, Y, \1, \w)$ be a vertex operator algebra. An anti-linear automorphism $\phi$ of $V$ is
an anti-linear isomorphism (as anti-linear map) $\phi:V\to V$ such
that $\phi(\bf1)=\bf1, \phi(\omega)=\omega$ and
$\phi(u_nv)=\phi(u)_n\phi(v)$ for any $ u, v\in V$ and $n\in
\mathbb{Z}$.}
\end{definition}
%Similarly, we could define the anti-linear isomorphism of the $V$-modules.
\begin{definition}{\rm
Let $(V, Y, \1, \w)$ be a vertex operator algebra and $\phi: V\to V$ be an anti-linear involution, i.e. an anti-linear automorphism of order $2$. The $(V, \phi)$ is called unitary if there exists a positive definite Hermitian form $(,): V\times V\to \mathbb{C}$
which is $\C$-linear on the first vector and anti-$\C$-linear on the second vector such that the following invariant property holds: for any $a, u, v\in V$\\
 $$(Y(e^{zL(1)}(-z^{-2})^{L(0)}a, z^{-1})u, v)=(u,
Y(\phi(a), z)v)$$where $L(n)$ is defined by $Y(\w, z)=\sum_{n\in \Z}L(n)z^{-n-2}$.}
\end{definition}
\begin{remark}
For a unitary vertex operator algebra $(V, \phi)$, the positive definite
Hermitian form $(,): V\times V\to \C$
is uniquely determined by the value $({\bf 1}, {\bf 1})$. In fact,
for any $u, v\in V_n$ there exists a complex number $\lambda\in \C$ such that
\begin{eqnarray*}
&&(u, v)=(u_{-1}{\bf 1}, v)\\
&&\ \ \ \ \ \ = (Y(u, z)\1, v)\\
&&\ \ \ \ \ \ = (\1, Y(e^{zL(1)}(-z^{-2})^{L(0)}\phi(u), z^{-1})v)\\
&&\ \ \ \ \ \  = (\1, Res_{z}z^{-1}Y(e^{zL(1)}(-z^{-2})^{L(0)}\phi(u), z^{-1})v)\\
&&\ \ \ \ \ \  = (\1, \lambda\1)\\
&&\ \ \ \ \ \ =\bar{\lambda}(\1,\1)
\end{eqnarray*}
We will
normalize the Hermitian form $(,)$ on $V$ such that $(\1, \1)=1$.
\end{remark}

\begin{definition}\label{dmodule}{\rm Let $(V, Y, \1, \w)$ be a vertex operator algebra and $\phi$ an anti-linear involution of $V,$ and  $g$ a finite order automorphism of $V$. An ordinary  $g$-twisted $V$-module $(M, Y_M)$ \cite{DLM2} is called a unitary
$g$-twisted $V$-module if there exists a positive definite
Hermitian form $(,)_M: M\times M\to \mathbb{C}$ which is $\C$-linear on the first vector and anti-$\C$-linear on the second vector such that the following invariant property:
$$(Y_M(e^{zL(1)}(-z^{-2})^{L(0)}a, z^{-1})w_1, w_2)_M=(w_1, Y_M(\phi(a),
z)w_2)_M$$ holds for $a \in V $ and $w_1, w_2\in M$.}
\end{definition}

%\begin{remark}
%Note that to satisfy the invariant property in the Definition $\ref{dmodule}$, $\phi$ should satisfy the property: $g\phi=\phi g^{-1}$.
%\end{remark}
Note that  if $(V, \phi)$ is a unitary vertex operator algebra, then $V$
is a unitary $V$-module.
\begin{lem}\label{unitary} Let $V$ be a vertex operator algebra and $\phi$ an anti-linear involution of $V,$ and  $g$ a finite order automorphism of $V.$ Then

(1) Any unitary $g$-twisted $V$-module $M$ is completely reducible.

(2) Any  unitary $g$-twisted $V$-module $M$  is a completely reducible module for the Virasoro algebra.
\end{lem}
\pf The proof of (1) is fairly standard using the invariant property. For (2) notice that the invariant property also implies that $(L(n)u, v)_M=(u, L(-n)v)_M$ for   $w_1, w_2\in M$ and $n\in \Z.$ As a result, $M$ is a unitary representation of the Virasoro algebra and the result follows immediately. \qed

\vskip0.5cm
In the following we construct unitary vertex operator algebras
from the given unitary vertex operator algebras. Recall that a vertex operator subalgebra $U=(U,Y,\1,\w')$
of $(V, Y, \1, \w)$ is a vector subspace $U$ of $V$ such that the restriction of $Y$ to $U$ gives a structure of vertex operator subalgebra on $U.$ The following proposition is immediate.
\begin{proposition}\label{unisub}
Let $(V, \phi)$ be a unitary vertex operator algebra and $U$ be a vertex operator subalgebra of $V$ such that the Virasoro element of $U$ is the same as that of $V$ and $\phi(U)=U$. Then $(U, \phi|_U)$ is a unitary vertex operator algebra.
\end{proposition}

Let $V$ be a vertex operator algebra and $g$ be a finite order automorphism of $V$. Then the fixed point subspace $V^g=\{a\in V| g(a)=a\}$ is a vertex operator subalgebra of $V.$
\begin{corollary}\label{obi}
Let $(V, \phi)$ be a unitary vertex operator algebra and $g$ be a finite
order automorphism of $V$ which commutes with $\phi$. Then
$(V^g, \phi|_{V^g})$ is a unitary vertex operator algebra.
\end{corollary}
\pf If $a$ lies in $V^G$ we
have $g\phi(a)=\phi g(a)=\phi(a).$ That is, $\phi(V^g)=V^g$. Then $(V^g, \phi|_{V^g})$ is unitary from Proposition \ref{unisub}.\qed

\vskip0.25cm
Let $(V, Y, \1, \w)$ be a vertex operator algebra and
$(U, Y, \1, \w')$ is a vertex operator subalgebra of $V$ such
that $\w'\in V_2$ and $L(1)\w'=0$. Then  $(U^c, Y, \1, \w-\w')$ is a vertex operator subalgebra of
$V$ where $U^c=\{v\in V|L'(-1)v=0\}$ \cite{FZ}.
\begin{corollary}
Let $(V, Y, \1, \w)$, $(U, Y, \1, \w')$ be vertex operator algebras satisfying the conditions above. Assume that $(V,\phi)$ is unitary and  $\phi(\w')=\w'$. Then $(U^c, \phi|_{U^c})$ is a unitary vertex operator algebra.
\end{corollary}
\pf  For $a\in U^c$, i.e.
$L'(-1)a=0$, we have
$L'(-1)\phi(a)=\w'_0\phi(a)=\phi(\w')_0\phi(a)=\phi(\w'_0a)=0.$ Thus $\phi(U^c)\subset U^c$. Since $\phi(\w')=\w'$, we have $\phi|_{U^c}$ is an anti-linear involution of $U^c$. Since
 $$(Y(e^{zL(1)}(-z^{-2})^{L(0)}a, z^{-1})u, v)=(u,
Y(\phi(a), z)v)$$
for $a, u, v\in U^c$, we have $$(Y(e^{zL''(1)}(-z^{-2})^{L''(0)}a, z^{-1})u, v)=(u,
Y(\phi(a), z)v)$$ holds for $a, u, v\in U^c$. Then $(U^c, \phi|_{U^c})$ is unitary.\qed

\vskip0.25cm
Now we recall some facts about the tensor product vertex operator algebra \cite{FHL}. Let $(V^1, Y_{V^1}, \1,
\omega^1),$ $ \cdots, (V^p, Y_{V^p}, \1, \omega^p)$ be vertex operator algebras. The tensor product of vertex operator algebras $V^1
, \cdots, V^p$ is constructed on
the tensor product vector space$$V = V^1\otimes \cdots \otimes
V^p$$ where the vertex operator $Y_{V}$ is defined by
$$Y_{V}(v^1\otimes \cdots \otimes v^p, z) = Y_{V^1}(v^1, z)\otimes
\cdots \otimes Y_{V^p}(v^p, z)$$\\for $v^i \in V^i\ (1\leq i \leq p),$
the vacuum vector is $$\1 = \1\otimes \cdots \otimes \1$$
and the Virasoro
element is
$$\omega = \omega^1 \otimes \cdots \otimes \1 + \cdots +\1\otimes
\cdots\otimes \omega^p.$$
Then  $(V, Y_{V}, \1, \omega)$
is a vertex operator algebra (see  \cite{FHL}, \cite{LL}).

Let $(M^i, Y_{M^i})$  be an ordinary $V^i$-module for $i=1,...,p.$  We may construct
the tensor product ordinary module $M^1\otimes \cdots \otimes M^p$ for the
tensor product vertex operator algebra $V^1\otimes \cdots \otimes
V^p$ by
$$Y_{M^1\otimes \cdots \otimes M^p}(v^1\otimes \cdots \otimes v^p, z) =
Y_{M^1}(v^1, z)\otimes \cdots \otimes Y_{M^p}(v^p, z).$$ Then  $(M^1\otimes \cdots \otimes M^p, Y_{M^1\otimes \cdots \otimes M^p})$ is an ordinary $V^1\otimes
\cdots \otimes V^p$-module.

\vskip0.25cm
Now let $(V^1, \phi_1),..., (V^p, \phi_p)$ be unitary vertex operator algebras and $(,)_i$ be the corresponding Hermitian form on $V^i$. We define a Hermitian form on $V^1\otimes\dots\otimes V^p$ as follow: for any $u_1\otimes\dots\otimes u_p, v_1\otimes\dots\otimes v_p \in  V^1\otimes\dots\otimes V^p$,$$(,): V^1\otimes\dots\otimes V^p \times V^1\otimes\dots\otimes V^p\to \mathbb{C}$$ $$(u_1\otimes\dots\otimes u_p, v_1\otimes\dots\otimes v_p)\mapsto(u_1, v_1)_1\dots (u_p, v_p)_p.$$ We also define an anti-linear map $\phi$ by: $$\phi: V^1\otimes\dots\otimes V^p\to V^1\otimes\dots\otimes V^p$$ $$\phi(u_1\otimes\dots\otimes u_p)\mapsto\phi_1(u_1)\otimes\dots\otimes \phi_p(u_p).$$ Obviously, (,) is a positive definite Hermitian form on $V^1\otimes\dots\otimes V^p $ and $\phi$ is an anti-linear involution of $V^1\otimes\dots\otimes V^p $. Now we have the following:
\begin{proposition}\label{tensor}
Let $(V^1, \phi_1),..., (V^p, \phi_p)$ be unitary vertex operator algebras and $\phi$ be the anti-linear involution of $V^1\otimes\dots\otimes V^p $ defined above. Then $(V^1\otimes\dots\otimes V^p, \phi)$ is a unitary vertex operator algebra.
\end{proposition}
\pf %Let $L^i(n)=\omega^i_{n+1}$, then we have
 %\begin{eqnarray*}
%&&(L(n)(u_1\otimes\dots\otimes u_p), v_1\otimes\dots\otimes v_p)\\
%&&=(L^1(n)u_1\otimes\dots \otimes u_p+\dots + u_1\otimes\dots\otimes L^p(n)u_p, v_1\otimes\dots\otimes v_p)\\
%&&=(L^1(n)u_1\otimes\dots \otimes u_p, v_1\otimes\dots\otimes v_p)+\dots +(u_1\otimes\dots\otimes L^p(n)u_p, v_1\otimes\dots\otimes v_p)\\
%&&=(L^1(n)u_1,v_1)\dots(u_p, v_p)+\dots +(u_1,v_1)\dots (L^p(n)u_p,v_p)\\
%&&=(u_1,L^1(-n)v_1)\dots(u_p, v_p)+\dots +(u_1,v_1)\dots (u_p,L^p(-n)v_p)\\
%&&=(u_1\otimes\dots \otimes u_p, L^1(-n)v_1\otimes\dots\otimes v_p)+\dots +(u_1\otimes\dots\otimes u_p, v_1\otimes\dots\otimes L^p(-n) v_p)\\
%&&=(u_1\otimes\dots\otimes u_p, L(n)(v_1\otimes\dots\otimes v_p)).
%\end{eqnarray*}
It is good enough to check that the invariant property
\begin{eqnarray*}
&&(Y_V(e^{zL(1)}(-z^{-2})^{L(0)}a_1\otimes\dots\otimes a_p, z^{-1})u_1\otimes\dots\otimes u_p, v_1\otimes\dots \otimes v_p)\\
&&=(u_1\otimes\dots\otimes u_p, Y_V(\phi(a_1\otimes\dots\otimes a_p), z)v_1\otimes\dots \otimes v_p)
\end{eqnarray*}
holds for $ a_1\otimes\dots\otimes a_p, u_1\otimes\dots\otimes u_p, v_1\otimes\dots \otimes v_p\in V=V^1\otimes\dots \otimes V^p$.

In fact, we have
\begin{eqnarray*}
&&(Y_V(e^{zL(1)}(-z^{-2})^{L(0)}a_1\otimes\dots\otimes a_p, z^{-1})u_1\otimes\dots\otimes u_p, v_1\otimes\dots \otimes v_p)\\
&&=(Y_{V^1}(e^{zL^1(1)}(-z^{-2})^{L^1(0)}a_1, z^{-1})u_1\otimes\dots \otimes\\
&&\ \ \ Y_{V^p}(e^{zL^p(1)}(-z^{-2})^{L^p(0)}a_p, z^{-1})u_p, v_1\otimes\dots \otimes v_p)\\
&&=(Y_{V^1}(e^{zL^1(1)}(-z^{-2})^{L^1(0)}a_1, z^{-1})u_1, v_1)_1\dots\\
&&\ \ \ (Y_{V^p}(e^{zL^p(1)}(-z^{-2})^{L^p(0)}a_p, z^{-1})u_p, v_p)_p\\
&&=(u_1, Y_{V^1}(\phi_1(a_1), z)v_1)_1\dots (u_p, Y_{V^p}(\phi_p(a_p), z)v_p)_p\\
&&=(u_1\otimes\dots\otimes u_p, Y_V(\phi(a_1\otimes\dots\otimes a_p), z)v_1\otimes\dots \otimes v_p).
\end{eqnarray*}
Then $(V^1\otimes\dots \otimes V^p, \phi)$ is a unitary vertex operator algebra.\qed

We could obtain the following proposition by the similar discussion as Proposition \ref{tensor}.
\begin{proposition}\label{tensm}
Let $V^1,..., V^p$ be vertex operator algebras and $\phi_i$ be an anti-linear involution of $V^i$  $(i=1,..., p)$. Assume that $M^i$ is a unitary module of $V^i\ (i=1,..., p)$, then $M^1\otimes \dots \otimes M^p$ is a unitary
$V^1\otimes\dots\otimes V^p$-module.
\end{proposition}
The following proposition is useful to prove the unitarity of vertex operator algebra.
\begin{proposition}\label{basic}
Let $V$ be a vertex operator algebra equipped with a positive
definite Hermitian form $(, ): V\times V\to \mathbb{C}$ and $\phi$ be an anti-linear involution of $V$. Assume that $V$ is
generated by the subset $S\subset V$, i.e.
$$V=span\{u^1_{n_1}\cdots u^k_{n_k}{\bf 1}| k\in \N, u^1,\cdots,
u^k\in S\}$$ and the invariant property
\begin{equation*} (Y(e^{zL(1)}(-z^{-2})^{L(0)}a, z^{-1})u,
v)=(u, Y(\phi(a), z)v)
 \end{equation*}
 holds for $a\in S, u, v\in V.$ Then $(V, \phi)$ is a unitary vertex operator algebra.
\end{proposition}
\pf Let $U$ be the subset of $V$ defined as follow:
$$U=\{a\in V|(Y(e^{zL(1)}(-z^{-2})^{L(0)}a, z^{-1})u, v)=(u, Y(\phi(a), z)v),\ \ \forall u,v \in V
\}.$$ It is easy to prove that ${\bf 1}\in U$. Now we
prove that if $a, b\in U$, then $a_nb\in U$ for any $n\in
\mathbb{Z}$. First, we have the following identity which was proved in  Theorem 5.2.1 of \cite{FHL},
\begin{eqnarray*}
&&-z_0^{-1}\delta(\frac{z_2-z_1}{-z_0})Y(e^{z_1L(1)}(-z_1^{-2})^{L(0)}a, z_1^{-1})Y(e^{z_2L(1)}(-z_2^{-2})^{L(0)}b, z_2^{-1})\\
&&+z_0^{-1}\delta(\frac{z_1-z_2}{z_0})Y(e^{z_2L(1)}(-z_2^{-2})^{L(0)}b, z_2^{-1})Y(e^{z_1L(1)}(-z_1^{-2})^{L(0)}a, z_1^{-1})\\
&&=z_1^{-1}\delta(\frac{z_2+z_0}{z_1})Y(e^{z_2L(1)}(-z_2^{-2})^{L(0)}Y(a, z_0)b, z_2^{-1}).
\end{eqnarray*}
By this identity, we have
\begin{eqnarray*}
&&Y(e^{z_2L(1)}(-z_2^{-2})^{L(0)}a_nb, z_2^{-1})\\
&&=Res_{z_1}\{-(-z_2+z_1)^nY(e^{z_1L(1)}(-z_1^{-2})^{L(0)}a, z_1^{-1})Y(e^{z_2L(1)}(-z_2^{-2})^{L(0)}b, z_2^{-1})\\
&&\ \ \ +(z_1-z_2)^nY(e^{z_2L(1)}(-z_2^{-2})^{L(0)}b, z_2^{-1})Y(e^{z_1L(1)}(-z_1^{-2})^{L(0)}a, z_1^{-1})\}.
\end{eqnarray*}
Then \begin{eqnarray*}
&&(Y(e^{z_2L(1)}(-z_2^{-2})^{L(0)}a_nb, z_2^{-1})u, v)\\
&&=Res_{z_1}\{(-(-z_2+z_1)^nY(e^{z_1L(1)}(-z_1^{-2})^{L(0)}a, z_1^{-1})Y(e^{z_2L(1)}(-z_2^{-2})^{L(0)}b, z_2^{-1})u, v)\\
&&\ \ \ +((z_1-z_2)^nY(e^{z_2L(1)}(-z_2^{-2})^{L(0)}b, z_2^{-1})Y(e^{z_1L(1)}(-z_1^{-2})^{L(0)}a, z_1^{-1})u,v)\}\\
&&=Res_{z_1}\{(u, -(-z_2+z_1)^nY(\phi(b), z_2)Y(\phi(a), z_1)v)\\
&&\ \ \ +(u, (z_1-z_2)^nY(\phi(a), z_1)Y(\phi(b), z_2)v)\}.
\end{eqnarray*}
On the other hand, we have
\begin{eqnarray*}
&&Y(\phi(a)_n\phi(b), z_2)\\
&&=Res_{z_1}\{(z_1-z_2)^nY(\phi(a), z_1)Y(\phi(b), z_2)-(-z_2+z_1)^nY(\phi(b), z_2)Y(\phi(a), z_1)\}.
\end{eqnarray*}
Thus we have
\begin{eqnarray*}
&&(Y(e^{z_2L(1)}(-z_2^{-2})^{L(0)}a_nb, z_2^{-1})u, v)\\
&&=(u, Y(\phi(a)_n\phi(b), z_2)v)\\
&&=(u, Y(\phi(a_nb), z_2)v).
\end{eqnarray*}
Then  we have $a_nb\in U$ for any $n\in
\mathbb{Z}$ if $a, b\in U$. Since $S\subset
U$ and $S$ generates $V$, we have $U=V$. So $(V, \phi)$ is a unitary vertex operator
algebra.\qed

\vskip0.25cm
The following proposition could be proved by the similar discussion as in Proposition \ref{basic}.
\begin{proposition}\label{module}
Let $g$ be a finite order automorphism of $V$ and  $M$ be an ordinary
$g$-twisted $V$-module equipped with a positive definite Hermitian
form $(, )_M: M\times M\to \mathbb{C}$. Assume that $V$ is a
vertex operator algebra which has an anti-linear involution $\phi$ and that $V$ is generated by the subset $S\subset
V$, i.e.
$$V=span\{u^1_{n_1}\cdots u^k_{n_k}{\bf 1}| k\in \N, u^1,\cdots,
u^k\in S\}$$ and the invariant property
\begin{equation*} (Y(e^{zL(1)}(-z^{-2})^{L(0)}a, z^{-1})w_1,
w_2)_M=(w_1, Y(\phi(a), z)w_2)_M
 \end{equation*}
 holds for $ a\in S,  w_1, w_2\in M$. Then $M$ is a unitary $g$-twisted $V$-module.
\end{proposition}
\section{Extension  vertex operator algebra by irreducible unitary simple current module}
\def\theequation{3.\arabic{equation}}
\setcounter{equation}{0} In this section, we assume that $(V, \phi)$ is
a unitary vertex operator algebra. For a unitary irreducible
$V$-module $M$, we construct an intertwining operator of type
$\left(\begin{tabular}{c}
$V$\\
$M$ $M$\\
\end{tabular}\right)$ by using the positive definite Hermitian forms $(,)_V,$ $(,)_M$ on $V, M$. We then use
intertwining operator to give a unitary vertex operator algebra structure on $V\oplus M$ if we further assume that
 $V$ is rational and $C_2$-cofinite, and $M$ is a simple current. These results will be useful in Section 4.

We first recall the  notion of intertwining operators  from [FHL].
\begin{definition}
Let $M^1$, $M^2$, $M^3$ be weak $V$-modules. An intertwining
operator $\mathcal {Y}( \cdot , z)$ of type $\left(\begin{tabular}{c}
$M^3$\\
$M^1$ $M^2$\\
\end{tabular}\right)$ is a linear map$$\mathcal
{Y}(\cdot, z): M^1\rightarrow Hom(M^2, M^3)\{z\}$$ $$v^1\mapsto
\mathcal {Y}(v^1, z) = \sum_{n\in \mathbb{C}}{v_n^1z^{-n-1}}$$
satisfying the following conditions:

(1) For any $v^1\in M^1, v^2\in M^2$and $\lambda \in \mathbb{C},
v_{n+\lambda}^1v^2 = 0$ for $n\in \mathbb{Z}$ sufficiently large.

(2) For any $a \in V, v^1\in M^1$,
$$z_0^{-1}\delta(\frac{z_1-z_2}{z_0})Y_{M^3}(a, z_1)\mathcal
{Y}(v^1, z_2)-z_0^{-1}\delta(\frac{z_2-z_1}{-z_0})\mathcal{Y}(v^1,
z_2)Y_{M^2}(a, z_1)$$
$$=z_2^{-1}\delta(\frac{z_1-z_0}{z_2})\mathcal{Y}(Y_{M^1}(a, z_0)v^1, z_2).$$

(3) For $v^1\in M^1$, $\dfrac{d}{dz}\mathcal{Y}(v^1,
z)=\mathcal{Y}(L(-1)v^1, z)$.
\end{definition}

In the following we assume that $M$ is a unitary irreducible $V$-module which has an anti-linear map $\psi$ such that $\psi(v_nw)=\phi(v)_n\psi(w)$ for $v\in V,$ $w\in M,$ $n\in\Z$  and $\psi^2=id$. Note that  $Y_M(., z)$ is an intertwining operator of type $\left(\begin{tabular}{c}
$M$\\
$V$ $M$\\
\end{tabular}\right)$. Define an operator $$\mathcal{Y}^*(., z):M\to Hom(V, M)\{ z\}$$ by the formula: for $ v\in V, w\in M$,

$$\mathcal{Y}^*(w, z)v=e^{zL(-1)}Y_M(v, -z)w.$$It is well-known that the operator $\mathcal{Y}^*(., z)$ is an intertwining operator of type $\left(\begin{tabular}{c}
$M$\\
$M$ $V$\\
\end{tabular}\right)$. We also define an operator $$\mathcal{Y}'(., z):M\to Hom(M,V)\{z\}$$ by the formula: for $ v\in V,  w_1, w_2\in M$,
$$(\mathcal{Y}'(w_1, z)w_2, v)_V=(w_2, \mathcal{Y}^*(e^{zL(1)}(-z^{-2})^{L(0)}\psi(w_1), z^{-1})v)_M.$$
\begin{proposition}
$\mathcal{Y}'(., z)$ is an intertwining operator of type $\left(\begin{tabular}{c}
$V$\\
$M$ $M$\\
\end{tabular}\right)$.
\end{proposition}
\pf Recall the following identities from
\cite{FHL}: for $f(z)\in z\mathbb{C}[[z]]$,

$$L(-1)e^{f(z)L(0)}=e^{f(z)L(0)}L(-1)e^{-f(z)},$$
$$L(1)e^{f(z)L(0)}=e^{f(z)L(0)}L(1)e^{f(z)},$$
\begin{eqnarray*}
&&L(-1)e^{f(z)L(1)}=e^{f(z)L(1)}L(-1)-2f(z)L(0)e^{f(z)L(1)}-f(z)^2L(1)e^{f(z)L(1)}\\
&&\ \ \ \ \ =e^{f(z)L(1)}L(-1)-2f(z)e^{f(z)L(1)}L(0)+f(z)^2e^{f(z)L(1)}L(1).
\end{eqnarray*}
{\bf Claim:} $\frac{d}{dz}\mathcal{Y}'(w_1, z)=\mathcal{Y}'(L(-1)w_1, z)$ for $w_1\in M$.

For any $v\in V, w_1, w_2\in M$, we have
\begin{eqnarray*}
&&(\frac{d}{dz}\mathcal{Y}'(w_1, z)w_2, v)_V=\frac{d}{dz}(w_2, \mathcal{Y}^*(e^{zL(1)}(-z^{-2})^{L(0)}\psi(w_1), z^{-1})v)_M\\
&&\ \ \ \ =(w_2, \mathcal{Y}^*(\frac{d}{dz}e^{zL(1)}(-z^{-2})^{L(0)}\psi(w_1), z^{-1})v)_M\\
&&\ \ \ \ \ +(w_2, \frac{d}{dz}\mathcal{Y}^*(w, z^{-1})|_{w=e^{zL(1)}(-z^{-2})^{L(0)}\psi(w_1)}v)_M.
\end{eqnarray*}

Since $\mathcal{Y}^*(., z)$ is an intertwining operator, we  prove the following identity by the similar discussion as in  Theorem 5.2.1 of \cite{FHL}:
\begin{eqnarray*}
&&\frac{d}{dz}\mathcal{Y}^*(w, z^{-1})|_{w=e^{zL(1)}(-z^{-2})^{L(0)}\psi(w_1)}=\mathcal{Y}^*(e^{zL(1)}(-z^{-2})^{L(0)}L(-1)\psi(w_1), z^{-1})\\
&&\ \ \ +\mathcal{Y}^*(2z^{-1}e^{zL(1)}L(0)(-z^{-2})^{L(0)}\psi(w_1), z^{-1})\\
&&\ \ \ -\mathcal{Y}^*(L(1)e^{zL(1)}(-z^{-2})^{L(0)}\psi(w_1), z^{-1}).
\end{eqnarray*}
Using the following identity (see \cite {FHL}):
  \begin{eqnarray*}
&&\frac{d}{dz}e^{zL(1)}(-z^{-2})^{L(0)}\\
&&=L(1)e^{zL(1)}(-z^{-2})^{L(0)}-2z^{-1}e^{zL(1)}L(0)(-z^{-2})^{L(0)},
\end{eqnarray*}
gives:
\begin{eqnarray*}
&&(\frac{d}{dz}\mathcal{Y}'(w_1, z)w_2, v)_V=(w_2, \mathcal{Y}^*(L(1)e^{zL(1)}(-z^{-2})^{L(0)}\psi(w_1), z^{-1})v)_M\\
&&\ \ \ \ \ -(w_2, \mathcal{Y}^*(2z^{-1}e^{zL(1)}L(0)(-z^{-2})^{L(0)}\psi(w_1), z^{-1})v)_M\\
&&\ \ \ \ \ +(w_2, \mathcal{Y}^*(e^{zL(1)}(-z^{-2})^{L(0)}L(-1)\psi(w_1), z^{-1})v)_M\\
&&\ \ \ \ \ +(w_2, \mathcal{Y}^*(2z^{-1}e^{zL(1)}L(0)(-z^{-2})^{L(0)}\psi(w_1), z^{-1})v)_M\\
&&\ \ \ \ \ -(w_2, \mathcal{Y}^*(L(1)e^{zL(1)}(-z^{-2})^{L(0)}\psi(w_1), z^{-1})v)_M\\
&&\ \ =(w_2, \mathcal{Y}^*(e^{zL(1)}(-z^{-2})^{L(0)}L(-1)\psi(w_1), z^{-1})v)_M\\
&&\ \ =(\mathcal{Y}'(L(-1)w_1, z)w_2, v)_V.
\end{eqnarray*}

{\bf Claim:} For any $v\in V$ and $w_1, w_2\in M$, we have
$$z_0^{-1}\delta(\frac{z_1-z_2}{z_0})Y(v, z_1)\mathcal {Y}'(w_1,
z_2)w_2-z_0^{-1}\delta(\frac{z_2-z_1}{-z_0})\mathcal{Y}'(w_1,
z_2)Y_M(v, z_1)w_2$$
$$=z_2^{-1}\delta(\frac{z_1-z_0}{z_2})\mathcal{Y}'(Y_{M}(v, z_0)w_1, z_2)w_2.$$
For any $v_1\in V$, we have the following identity which was
essentially proved in  Theorem 5.2.1 of \cite{FHL},
\begin{eqnarray*}
&&-z_0^{-1}\delta(\frac{z_2-z_1}{-z_0})Y_M(e^{z_1L(1)}(-z_1^{-2})^{L(0)}v, z_1^{-1})\mathcal{Y}^*(e^{z_2L(1)}(-z_2^{-2})^{L(0)}w_1, z_2^{-1})v_1\\
&&+z_0^{-1}\delta(\frac{z_1-z_2}{z_0})\mathcal{Y}^*(e^{zL(1)}(-z_2^{-2})^{L(0)}w_1, z_2^{-1})Y(e^{z_1L(1)}(-z_1^{-2})^{L(0)}v, z_1^{-1})v_1\\
&&=z_2^{-1}\delta(\frac{z_1-z_0}{z_2})\mathcal{Y}^*(e^{z_2L(1)}(-z_2^{-2})^{L(0)}Y_M(v, z_0)w_1, z_2^{-1})v_1.
\end{eqnarray*}
Since $M$ is a unitary $V$-module, we have
\begin{eqnarray*}
&&(z_0^{-1}\delta(\frac{z_1-z_2}{z_0})Y(v, z_1)\mathcal
{Y}'(w_1, z_2)w_2, v_1)_V\\
&&\ \ \ \ -(z_0^{-1}\delta(\frac{z_2-z_1}{-z_0})\mathcal{Y}'(w_1,
z_2)Y_{M}(v, z_1)w_2, v_1)_V\\
&&=(w_2,z_0^{-1}\delta(\frac{z_1-z_2}{z_0})\mathcal{Y}^*(e^{z_2L(1)}(-z_2^{-2})^{L(0)}\psi(w_1), z_2^{-1})\cdot\\
&&\ \ \ \ \ \cdot Y(e^{z_1L(1)}(-z_1^{-2})^{L(0)}\phi(v), z_1^{-1})v_1 )_M\\
&&\ \ \ -(w_2, z_0^{-1}\delta(\frac{z_2-z_1}{-z_0})Y_M(e^{z_1L(1)}(-z_1^{-2})^{L(0)}\phi(v), z_1^{-1})\cdot\\
&&\ \ \ \ \ \cdot \mathcal{Y}^*(e^{z_2L(1)}(-z_2^{-2})^{L(0)}\psi(w_1), z_2^{-1})v_1)_M\\
&&=(w_2,z_2^{-1}\delta(\frac{z_1-z_0}{z_2})\mathcal{Y}^*(e^{z_2L(1)}(-z_2^{-2})^{L(0)}Y_M(\phi(v), z_0)\psi(w_1), z_2^{-1})v_1 )_M\\
&&=(w_2,z_2^{-1}\delta(\frac{z_1-z_0}{z_2})\mathcal{Y}^*(e^{z_2L(1)}(-z_2^{-2})^{L(0)}\psi(Y_M(v, z_0)w_1), z_2^{-1})v_1 )_M\\
&&=(z_2^{-1}\delta(\frac{z_1-z_0}{z_2})\mathcal{Y}'(Y_M(v, z_0)w_1,z_2),v_1)_V.
\end{eqnarray*}
This completes the proof. \qed

\vskip0.5cm
We now assume that $(V, \phi)$ is a rational and
$C_2$-cofinite unitary vertex operator algebra. Recall that
an irreducible $V$-module $M$ is called simple current if the
tensor product $M\boxtimes M_1$ is an irreducible $V$-module for
any irreducible $V$-module $M_1$. It was proved in \cite{LY} that there exists a vertex operator algebra structure on $U=V\oplus M$ if $M$ is a simple current $V$-module satisfying some additional conditions. By using the intertwining operator constructed above, we construct a unitary vertex operator algebra structure on $U=V\oplus M$ in the following way.
\begin{theorem}\label{extension}
Let $(V, \phi)$ be a rational and $C_2$-cofinite unitary self-dual vertex
operator algebra and $M$ be a simple current
irreducible $V$-module having integral weights. Assume that $M$
has an anti-linear map $\psi$ such that
$\psi(v_nw)=\phi(v)_n\psi(w)$ and $\psi^2=id$,
$(\psi(w_1), \psi(w_2))_M=\overline{(w_1, w_2)_M}$ and the
Hermitian form $(,)_V$  on $V$ has the property that $(\phi(v_1),
\phi(v_2))_V=\overline{(v_1, v_2)_V}$. Then
$(U, \phi_U)$ has a unique unitary vertex operator algebra structure, where $\phi_U:U\to U$ is the anti-linear involution defined by $\phi_U(v, w)=(\phi(v), \psi(w))$, for $ v\in V,  w\in M$. Furthermore, $U$ is rational and $C_2$-cofinite.
\end{theorem}
\pf %%First, let $\phi_U$ be the anti-linear map $\phi_U:U\to U$ by $\phi_U(v, w)=(\phi(v), \psi(w))$, for $\forall v\in V, \forall w\in M$.
If $U$ has a vertex operator algebra structure, then the vertex operator algebra structure is unique \cite{DM2}. Moreover, $U$ is rational and $C_2$-cofinite \cite{Y}. So it is good enough to construct a unitary vertex operator algebra structure on $U$. Let $(,)_U: U\times U\to \mathbb{C}$ be the Hermitian form on $U$ defined by $(v_1, v_1)_U=(v_1, v_2)_V$, $(v_1, w_1)=0$ and $(w_1, w_2)_U=(w_1, w_2)_M$ for any $v_1, v_2\in V$, $w_1, w_2\in M.$ It is obvious that this Hermitian form is positive definite.

Define a linear operator $Y_U(., z): U\to End(U)[[z^{-1}, z]]$ by $$Y_U(v, z)=\left(\begin{tabular}{cccc}
$Y(v, z)$ $ $ $ $$0$\\
$0$ $ $ $ $$Y_M(v, z)$\\
\end{tabular}\right)$$

$$Y_U(w, z)=\left(\begin{tabular}{cccc}
0 $ $ $ $$\mathcal{Y}'(w, z)$\\
$\mathcal{Y}^*(w, z)$ $ $ $ $ 0\\
\end{tabular}\right)$$
for any $v\in V$ and $w\in M$.

{\bf Claim:} The operator $Y_U(., z)$ satisfies the skew-symmetry property, i.e. for any $u,v\in U$, we have $$Y_U(u, z)v=e^{zL(-1)}Y_U(v, -z)u.$$
We need the following identities:
$$(-z^2)^{L(0)}e^{zL(1)}(-z^2)^{-L(0)}=e^{-z^{-1}L(1)},$$
$$(-z^2)^{-L(0)}e^{zL(-1)}(-z^2)^{L(0)}=e^{-z^{-1}L(-1)}.$$
Here is a proof of the second identity and the proof of the first identity is similar. It is enough to show
 $$(-z^2)^{-L(0)}zL(-1)(-z^2)^{L(0)}=-z^{-1}L(-1)$$ or
 $$z^{-L(0)}L(-1)z^{L(0)}=z^{-1}L(-1)$$
which is clear.

Now we prove the claim.  By definition we need to show that
$$Y_U(w_1, z)w_2=e^{zL(-1)}Y_U(w_2, -z)w_1$$ for $w_1, w_2\in M$.
For any $v_1\in V$, we have
\begin{eqnarray*}
&&(Y_U(w_1, z)w_2, v_1)_V\\
&&=(w_2, \mathcal{Y}^*(e^{zL(1)}(-z^{-2})^{L(0)}\psi(w_1), z^{-1})v_1)_M\\
&&=(w_2, e^{z^{-1}L(-1)}Y_M(v_1, -z^{-1})e^{zL(1)}(-z^{-2})^{L(0)}\psi(w_1))_M
\end{eqnarray*}
and
\begin{eqnarray*}
&&(e^{zL(-1)}Y_U(w_2, -z)w_1, v_1)_V\\
&&=(w_1, \mathcal{Y}^*(e^{-zL(1)}(-z^{-2})^{L(0)}\psi(w_2), -z^{-1})e^{zL(1)}v_1)_M\\
&&=(w_1, e^{(-z^{-1}L(-1))}Y_M(e^{zL(1)}v_1, z^{-1})e^{-zL(1)}(-z^{-2})^{L(0)}\psi(w_2))_M\\
&&=(Y_M(e^{z^{-1}L(1)}(-z^2)^{L(0)}\phi(e^{zL(1)}v_1), z)e^{-z^{-1}L(1)}w_1,e^{-zL(1)}(-z^{-2})^{L(0)}\psi(w_2) )_M\\
&&=(Y_M((-z^2)^{L(0)}\phi(v_1), z)e^{-z^{-1}L(1)}w_1,e^{-zL(1)}(-z^{-2})^{L(0)}\psi(w_2) )_M\\
&&=((-z^{-2})^{L(0)}e^{-zL(-1)}Y_M((-z^2)^{L(0)}\phi(v_1), z)e^{-z^{-1}L(1)}w_1,\psi(w_2) )_M\\
&&=(e^{z^{-1}L(-1)}(-z^{2})^{-L(0)}Y_M((-z^2)^{L(0)}\phi(v_1), z)e^{-z^{-1}L(1)}w_1,\psi(w_2) )_M\\
&&=(e^{z^{-1}L(-1)}Y_M(\phi(v_1), -z^{-1})(-z^{2})^{-L(0)}e^{-z^{-1}L(1)}w_1,\psi(w_2) )_M\\
&&=(e^{z^{-1}L(-1)}Y_M(\phi(v_1), -z^{-1})e^{zL(1)}(-z^{2})^{-L(0)}w_1,\psi(w_2) )_M\\
&&=\overline{(\psi(e^{z^{-1}L(-1)}Y_M(\phi(v_1), -z^{-1})e^{zL(1)}(-z^{2})^{-L(0)}w_1),w_2 )_M}\\
&&=\overline{(e^{z^{-1}L(-1)}Y_M(v_1, -z^{-1})e^{zL(1)}(-z^{2})^{-L(0)}\psi(w_1),w_2 )_M}\\
&&=(w_2, e^{z^{-1}L(-1)}Y_M(v_1, -z^{-1})e^{zL(1)}(-z^{2})^{-L(0)}\psi(w_1))_M\\
&&=(Y_U(w_1, z)w_2, v_1)_V.
\end{eqnarray*}
So the claim is established.

We can now prove that $(U, Y_U(., z))$ is a vertex operator algebra. Since the skew symmetry holds, it is
enough to have the locality. By Theorem 5.6.2 of \cite{FHL}, we only need to prove the locality
for three elements in $M.$  But this follows a similar discussion in  Proposition 3
of \cite{LY}.

To prove that $(U, \phi_U)$ is a unitary vertex operator algebra, we first prove that $\phi_U$ is an anti-linear involution of $U$. Obviously, the order of $\phi_U$ is 2. So it suffices to prove that $\phi_U(u_nv)=\phi_U(u)_n\phi_U(v)$ for any $u, v\in U$. We now verify this property  case by case. If $u, v\in V$ this is obvious. If $u\in V$ and $v\in M$, we have $\phi_U(u_nv)=\phi(u)_n\psi(v)=\phi_U(u)_n\phi_U(v).$
  If $u\in M, v\in V$, we have $\phi_U(Y_U(u, z)v)=\phi_U(e^{zL(-1)}Y(v, -z)u)=e^{zL(-1)}Y(\phi(v), -z)\psi(u)=Y_U(\phi_U(u), z)\phi_U(v)$,
   this implies  $\phi_U(u_nv)=\phi_U(u)_n\phi_U(v).$
   If $u\in M, v\in M$, for any $ v_1\in V$ we have:
\begin{eqnarray*}
&&(\phi_U(Y_U(u, z)v), \phi(v_1))\\
&&=\overline{(Y_U(u, z)v, v_1)}\\
&&=\overline{(v, \mathcal{Y}^*(e^{zL(0)}(-z^{-2})^{L(0)}\psi(u), z^{-1})v_1)}\\
&&=(\psi(v), \psi(\mathcal{Y}^*(e^{zL(0)}(-z^{-2})^{L(0)}\psi(u), z^{-1})v_1))\\
&&=(\psi(v), \mathcal{Y}^*(e^{zL(0)}(-z^{-2})^{L(0)}u, z^{-1})\phi(v_1))\\
&&=(Y_U(\psi(u), z)\psi(v), \phi(v_1)).
\end{eqnarray*}
 Thus we have
$\phi_U(u_nv)=\phi_U(u)_n\phi_U(v)$, then $\phi_U$ is an
anti-linear involution of $U$.

It remains to prove the invariant property
$$(Y(e^{zL(1)}(-z^{-2})^{L(0)}a, z^{-1})u, v)=(u, Y(\phi_U(x),
z)v)$$ holds for any  $ a, u, v\in U$. This is obvious from the
definition of $Y_U(., z)$. The proof of theorem is complete.\qed
%\begin{remark}
%Note that in the Proposition 3 of \cite{LY}, the simple current
%module $M$ is assumed to be self-dual. Let $L$ is a positive
%definite even lattice, we'll prove that the vertex operator
%algebra $V_L^+$ is a unitary vertex operator algebra such that the
%anti-linear isomorphism satisfying the conditions in Theorem
%\ref{extension}, and we know there exists unitary irreducible
%$V_L^+$-module which is not self-dual (see \cite{ADL}). So Theorem
%\ref{extension} complements the Proposition 3 of \cite{LY}.
%\end{remark}
\section{Examples of unitary vertex operator algebras}
\def\theequation{4.\arabic{equation}}
\setcounter{equation}{0}

In this section we prove that most of the well-known rational
and $C_2$-cofinite vertex operator algebras are unitary vertex
operator algebras. However, there exist some unitary vertex
operator algebras which are neither rational nor $C_2$-cofinite.
\subsection{Unitary Virasoro vertex operator algebras}

In this subsection we construct unitary vertex operator algebras associated to the Virasoro algebra. First, we recall some facts about
Virasoro vertex operator algebras \cite{FZ}, \cite{W}. We denote the Virasoro algebra by  $L=\oplus_{n\in \mathbb{Z}}
\mathbb{C}L_n\oplus \mathbb{C}C$ with the
commutation relations $$[L_m, L_n]=(m-n)L_{m+n}+\frac{1}{12}(m^3-m)\delta_{m+n,
0}C,$$
$$[L_m, C]=0.$$

Set $\mathfrak{b}=(\oplus_{n\geq
1}\mathbb{C}L_n)\oplus(\mathbb{C}L_0\oplus \mathbb{C}C)$, then we
know that $\mathfrak{b}$ is a subalgebra of $L$. For any two
complex numbers $(c, h)\in \C$, let $\C$ be a 1-dimensional
$\mathfrak{b}$-module defined as follows:
$$L_n\cdot 1=0, n\geq 1,$$
$$L_0\cdot 1=h\cdot 1,$$
$$C\cdot 1=c\cdot 1.$$

Set $$V(c, h)=U(L)\otimes_{U(\mathfrak{b})}\C$$ where
$U(.)$ denotes the universal enveloping algebra. Then $V(c, h)$ is
a highest weight module of the Virasoro algebra of highest weight $(c, h)$, which is called the Verma module of Virasoro algebra,
and $V(c, h)$ has a unique maximal proper submodule $J(c, h)$. Let
$L(c, h)$ be the unique irreducible quotient module of $V(c, h)$.
Set
$$\overline{V(c, 0)}=V(c, 0)/(U(L)L_{-1}1\otimes 1),$$ it is well-known
that $\overline{V(c, 0)}$ has a vertex operator algebra structure with Virasoro element $\w=L_{-2}1$
and $ L(c, 0)$ is the unique irreducible quotient vertex operator
algebra of $\overline{V(c, 0)}$ \cite{FZ}.

For $m\geq 2$, set $$c_m=1-\frac{6}{m(m+1)},$$ $$h^m_{r,
s}=\frac{(r(m+1)-sm)^2-1}{4m(m+1)}, (1\leq s\leq r \leq m-1).$$
It was proved in \cite{W}, \cite{DLM3} that $L(c_m, 0)$ $(m\geq 2)$ are rational and
$C_2$-cofinite vertex operator algebra and $L(c_m, h^m_{r, s})$ are
the complete list of irreducible $L(c_m, 0)$-modules.

Now we recall some fact about the Hermitian form
on $L(c, h)$. For $(c, h)\in \R$, it was proved in Proposition 3.4
of \cite{KR} that there is a unique Hermitian form $(,)$
 such that $$(v_{c, h}, v_{c, h} )=1,$$$$(L_nu, v)=(u, L_{-n}v),$$
for any $u, v\in L(c, h)$, where $v_{c, h}$ denotes the highest weight vector of $L(c, h)$.
 It is well-known that $L(c, h)$ is unitary, i.e. the Hermitian
 form $(,)$
 on $L(c, h)$ is positive definite, if and only if $c\geq1, h\geq0$ or
 $c=c_m, h=h^m_{r,s}$ \cite{KR}.

\vskip0.25cm
For any real number $c$, define an anti-linear map $\overline{\phi}$ of
$\overline{V(c, 0)}$ as follow:

$$\overline{\phi}:\overline{V(c, 0)}\to \overline{V(c,
0)}$$
$$L_{-n_1}\cdots L_{-n_k}\cdot 1\mapsto L_{-n_1}\cdots
L_{-n_k}\cdot 1, n_1\geq \cdots \geq n_k\geq 2.$$

\begin{lem}\label{involution}
Assume that $c\in \R$, and let $\overline{\phi}$ be the anti-linear map defined above. Then $\overline{\phi}$ is an anti-linear involution of vertex operator algebra $\overline{V(c, 0)}.$ Furthermore, $\overline{\phi}$ induces an anti-linear involution $\phi$ of $L(c, 0).$
\end{lem}
\pf Since we have $\overline{\phi}^2=id$, it is good enough to prove that $\overline{\phi}$ is an anti-linear automorphism. Let $U$ be the subspace of $\overline{V(c, 0)}$ which is defined by $$U=\{u\in \overline{V(c, 0)}|\overline{\phi}(u_nv)=\overline{\phi}(u)_n\overline{\phi}(v), \forall v\in \overline{V(c, 0)}, n\in \Z\}.$$It is easy to prove that if $a, b\in U$, then $a_mb\in U$ for any $m\in \Z$. Note that $\1\in U$ and $\w=L_{-2}\cdot 1\in U$. Thus $U=\overline{V(c, 0)}$ as $\overline{V(c, 0)}$ is generated by $\w$. This implies that $\overline{\phi}$ is an anti-linear involution of $\overline{V(c, 0)}$.

Let $\overline{J(c, 0)}$ be the maximal
proper $L$-submodule of $\overline{V(c, 0)}$, we have
$\overline{\phi}(\overline{J(c, 0)})$ is a proper $L$-submodule of
$\overline{V(c, 0)}$, this implies
$\overline{\phi}(\overline{J(c, 0))}\subset\overline{J(c,
0)}$. Thus $\overline{\phi}$ induces an anti-linear involution
$\phi$ of $L(c, 0)$.\qed

Now we have the main result in the subsection.

\begin{theorem}\label{Vi1}
Assume that $c\in \R$ and let $\phi$ be the anti-linear involution of $L(c_m, 0)$ defined above.  Then $(L(c, 0),\phi)$  is a  unitary vertex operator algebra if and only if $c\geq1$ or $c=c_m$ for some integer $m\geq 2$.
\end{theorem}
\pf  Assume that $c\geq1$ or $c=c_m$ for some integer $m\geq 2$. Then the Hermitian form $(,)$ defined above is positive definite. We now prove the invariant property in the definition of the unitary vertex operator algebra. Recall from \cite{FZ}, \cite{W} that $L(n)u=L_nu$ for any $u\in L(c, 0)$, then we have $(L(n)u, v)=(u, L(-n)v)$ for any $u, v\in L(c, 0)$. Thus we have
\begin{eqnarray*}
&&(Y(e^{zL(1)}(-z^{-2})^{L(0)}\omega, z^{-1})u, v)=z^{-4}(Y(\omega, z^{-1})u, v)\\
&&\ \ \ \ =\sum_{n\in \Z}(\omega_{n+1}u, v)z^{n-2}\\
&&\ \ \ \ =\sum_{n\in \Z}(L(n)u, v)z^{n-2}\\
&&\ \ \ \ =\sum_{n\in \Z}(u, L(-n)v)z^{n-2}\\
&&\ \ \ \ =\sum_{n\in \Z}(u, \w_{-n+1}v)z^{n-2}\\
&&\ \ \ \ =(u, Y(\w, z)v)\\
&&\ \ \ \ =(u, Y(\phi(\w), z)v).
\end{eqnarray*}
Since $L(c, 0)$ is generated by $\omega$,  $(L(c, 0), \phi)$ is
a unitary vertex operator algebra by Proposition \ref{basic}.
%The similar calculation could prove that $L(c_m, h_{r, s})$ are the
%complete list of unitary irreducible $L(c_m, 0)$-modules.

We now assume that $(L(c, 0), \phi)$ is a unitary vertex operator algebra. By the Lemma \ref{unitary}, $L(c, 0)$ is unitary as the module of the Virasoro algebra, thus $c\geq1$ or $c=c_m$ for some integer $m\geq2$.
\qed
\begin{remark}
 Note that if $c\geq1,$ $(L(c, 0),\phi)$ is a unitary vertex operator algebra, although $L(c, 0)$  is neither rational nor $C_2$-cofinite.
 \end{remark}
We also have the following proposition by the similar discussion as in Theorem \ref{Vi1}.
\begin{theorem}
Assume that $c\in \R$ and let $\phi$ be the anti-linear involution of $L(c, 0)$ defined above. Then $L(c, h)$ is a unitary module of $L(c, 0)$ if and only if $c\geq 1, h\geq 0$ or $c=c_m, h=h^m_{r,s}$.
\end{theorem}
\subsection{Unitary affine vertex operator algebras}
In this subsection we construct unitary vertex operator algebra associated to the affine Kac-Moody
algebras. First,
we recall some facts about affine vertex operator algebras \cite{FZ}. Let
$\g$ be a finite dimensional simple Lie algebra and $\h$ a Cartan subalgebra. Fix
a non-degenerate symmetric invariant bilinear form $(,)$ on $\g$ so that $(\theta, \theta)=2$ where $\theta$ is the
maximal root of $\g$. Consider the affine Lie algebra
$\hat{\g}=\g\otimes\C[t, t^{-1}]\oplus \C K$ with the commutation relations
$$[x(m), y(n)]=[x, y](m+n)+Km(x, y)\delta_{m+n, 0},$$
$$[\hat{\g}, K]=0,$$ where $x(m)=x\otimes t^{n}$.

For complex number $k\in \C$, set
$$V_\g(k)=U(\hat{\g})/J_k,$$ where $J_k$ is the left ideal of
$U(\hat{\g})$ generated by $x(n)$ and $K-k$ for $x\in \g$ and $n\geq 0.$ It is
well-known that  $V_\g(k)$ has a vertex operator algebra structure if
$k\neq -h^{\vee}$ where $h^{\vee}$ is the dual Coxeter number of
$\g$. Moreover, $V_\g(k)$ has a unique maximal proper $\hat{\g}$-submodule
$J(k)$ and  $L_\g(k, 0)=V_\g(k)/J(k)$ is a simple vertex operator algebra. As usual
we denote the corresponding irreducible highest weight module for $\hat{\g}$ associated to a highest weight $\lambda\in\h^*$
of $\g$ by $L_\g(k, \lambda)$. It was proved in \cite{FZ}, \cite{DLM3} that $L_\g(k, 0)$ is a
simple rational and $C_2$-cofinite vertex operator algebra if
$k\in \Z^{+}$ and
$$\{L_\g(k, \lambda)|(\lambda, \theta)\leq k, \lambda\in \h^*\ is\ integral\ dominant\}$$ are the complete list of inequivalent irreducible $L_\g(k, 0)$-modules.
\vskip0.25cm
 Let $\omega_0$ be  the compact involution \cite{K} of  $\g$ which is the anti-linear automorphism determined by:
$$\omega_0(e_i)=-f_i,\ \ \omega_0(f_i)= -e_i,\ \ \omega_0( h_i)= -h_i,$$ where $\{h_i, e_i, f_i\}$ are the
Chevalley generators of $\g$.
\begin{lem}\label{kill}
For any $x, y\in \g,$ we have $\overline{(\w_0(x), \w_0(y))}=(x, y)$.
\end{lem}
\pf Note that it is good enough to prove that $\overline{\kappa(\w_0(x), \w_0(y))}=\kappa(x, y)$, where $\kappa(,)$ is the Killing form of $\g$. Recall that $\kappa(x, y)=tr(ad xad y)$, let $z_1, \cdots, z_k$ be the bases of $\g$ such that $adxad y$ is a upper-triangular matrix, then we have $adxady (z_i)=\lambda_iz_i+w_i$ for $1\leq i\leq k,$ where $w_i$ is some linearly combine of $z_1,\cdots, z_{i-1}$. On the other hand, we have $[\w_0(x), [\w_0(y), \w_0(z)]]=\w_0([x, [y, z]])=\w_0(adxad y(z))=\overline{\lambda_i}\w_0(z_i)+v_i$ for $1\leq i\leq k$, where $v_i$ is some linearly combine of $\w_0(z_1),\cdots, \w_0(z_{i-1})$. This implies that $\overline{\kappa(\w_0(x), \w_0(y))}=\kappa(x, y)$.
\qed

\vskip0.25cm
Let $\hat{\w}_0$ be the compact involution \cite{K}  of $\hat{\g}$ which is the anti-linear automorphism determined by:
$$\hat{\omega}_0(x\otimes t^{m})= \omega_0(x)\otimes t^{-m}, \hat{\omega}_0(K)=-K.$$
It was proved in Theorem 11.7 of \cite{K} that if $k\in \Z^{+}$ and $\lambda$ is an integral dominant weight such that $(\lambda, \theta)\leq k$, then $L_\g(k, \lambda) $ has a unique positive definite Hermitian form $(,)$ such
that $$({ 1}, { 1})=1,$$ $$(xu, v)=-(u,\hat{\omega}_0(x)
v)$$
for  $x\in \hat{\g},$  $u, v \in L_\g(k, \lambda)$ where $1$ is a fixed highest weight vector of $L_\g(k, \lambda).$

We also introduce a linear map $\w_{\hat{\g}}: \hat{\g}\to \hat{\g}$ of $\hat{\g}$ such that  $\w_{\hat{\g}}(x(n))=\w_0(x)(n)$ and $\w_{\hat{\g}}(K)=K.$
Let $T(\hat{\g}) $ be  the tensor product algebra of the affine Lie algebra $\hat{\g}$. Define an anti-linear map $\Phi_T$ of $T(\hat{\g})$ as follow: for any $y_1, ..., y_k\in \hat{\g}$$$\Phi_T:T(\hat{\g})\to T(\hat{\g})$$
$$y_1\otimes\cdots\otimes y_k\mapsto \w_{\hat{\g}}(y_1)\otimes\cdots\otimes \w_{\hat{\g}}(y_k).$$ It is easy to prove that $\Phi_T$ is an anti-linear involution such that $\Phi_Tx(m)\Phi_T^{-1}=\w_0(x)(m)$. Recall that the universal enveloping algebra $U(\hat{\g})$ is defined to be the quotient $T(\hat{\g})/I$, where $I$ is the ideal of $T(\hat{\g})$ generated by $a\otimes b-b\otimes a-[a, b]$, $ a, b\in \hat{\g}$. By Lemma \ref{kill}, we have
\begin{eqnarray*}
&&\Phi_T(x(m)\otimes y(n)-y(n)\otimes x(m)-[x(m), y(n)])\\
&&=\w_0(x)(m)\otimes \w_0(y)(n)-\w_0(y)(n)\otimes \w_0(x)(m)\\
&&\ \ \ -[\w_0(x), \w_0(y)](m+n)-Km\overline{(x, y)}\delta_{m+n, 0}\\
&&=\w_0(x)(m)\otimes \w_0(y)(n)-\w_0(y)(n)\otimes \w_0(x)(m)\\
&&\ \ \ -[\w_0(x), \w_0(y)](m+n)-Km(\w_0(x), \w_0(y))\delta_{m+n, 0}\\
&&=\w_0(x)(m)\otimes \w_0(y)(n)-\w_0(y)(n)\otimes \w_0(x)(m)-[\w_0(x)(m), \w_0(y)(n)].
\end{eqnarray*}
This implies that $\Phi_T(I)\subset I$, then $\Phi_T$ induces an anti-linear involution  $\Phi_U$ of $U(\hat{\g})$ such that $\Phi_Ux(m)\Phi_U^{-1}=\w_0(x)(m)$. Note that if $k\in \R$, we have $\Phi_U(J_k)\subset J_k$, then $\Phi_U$ induces an anti-linear map $\Phi$ of $V_\g(k)$ such that  $\Phi x(m)\Phi^{-1}=\w_0(x)(m)$ and $\Phi^2=id$.
\begin{lem}\label{involution1}
Assume that $k\in \R$ and $k\neq h^{\vee}$, and let $\Phi$ be the anti-linear map defined above. Then $\Phi$ is an anti-linear involution of vertex operator algebra $V_\g(k).$ Furthermore, $\Phi$ induces an anti-linear involution $\phi$ of $L_\g(k, 0).$
\end{lem}
\pf Since we have $\Phi^2=id$, it is good enough to prove that $\Phi$ is an anti-linear automorphism. Let $U$ be the subspace of $V_\g(k)$ which is defined by $$U=\{u\in V_\g(k)|\Phi(u_nv)=\Phi(u)_n\Phi(v), \forall v\in V_\g(k), n\in \Z\}.$$It is easy to prove that if $u, v\in U$, then $u_mv\in U$ for any $m\in \Z$. Note that $\1\in U$ and $x(-1)\1\in U$ for any $x\in \g$. Thus we have $U=V_\g(k)$, since $V_\g(k)$ is generated by $x(-1)\1, x\in \g$. Thus  $\Phi$ is an anti-linear involution of $V_\g(k)$.

Note that $\Phi(J(k))$ is a proper $\hat{\g}$-submodule of $V_\g(k).$ This forces $\Phi(J(k))\subset J(k)$. So
$\Phi$ induces an anti-linear involution $\phi$ of $L_g(k, 0)$.\qed

Now we can prove the main result in this subsection.
\begin{theorem}\label{aff1}
Assume that $k\in \R$ and $k\neq h^{\vee}$, and let $\phi$ be the anti-linear automorphism of $L_g(k, 0)$ defined as above. Then $(L_g(k, 0), \phi)$ is a unitary vertex operator algebra if and only if $k\in \Z^{+}$.
\end{theorem}
\pf Assume that $k\in \Z^{+}.$  Then the Hermitian form $(,)$ on $L_\g(k, 0)$ is positive definite. We now prove the invariant property. Recall from \cite{FZ} that $Y(x(-1)\1, z)=\sum_{n\in \Z}x(n)z^{-n-1}$ for $x\in\g.$ Note that  $$(x(n) u, v)=-(u,\w_0(x)(-n)v)$$ for $x\in \hat{\g}$ and $ u, v \in L_\g(k, 0).$
Thus we have:
\begin{eqnarray*}
&&(Y(e^{zL(1)}(-z^{-2})^{L(0)}x(-1)\1, z^{-1})u, v)\\
&&=(Y(-z^{-2}x(-1)\1, z^{-1})u, v)\\
&&=\sum_{n\in \Z}-z^{-2}(x(n)u, v)z^{n+1}\\
&&=\sum_{n\in \Z}(u, \w_0(x)(-n)v)z^{n-1}\\
&&=(u, Y(\phi(x(-1)\1), z)v).
\end{eqnarray*}
Since $L_\g(k, 0)$ is generated by $x(-1)\1$, we have $(L_\g(k, 0),\phi)$ is a unitary vertex operator algebra by Proposition \ref{basic}.

Conversely, assume that $(L_\g(k, 0), \phi)$ is a unitary vertex operator algebra. By the definition and the calculation as above, we have a positive definite Hermitian form $(,)$ on $L_\g(k, 0)$ such that for any $x\in \g$ $$(x(n)u, v)=-(u, \w_0(x)(-n)v).$$ This implies that $L_\g(k, 0)$ is unitary as $\hat{\g}$-module and $k\in \Z^{+}$ by Theorem 11.7 of \cite{K}.\qed

The following result is immediate  by the similar discussion as in Theorem \ref{aff1}.
\begin{theorem}
Assume that $k\in \Z^{+}$  and let $\phi$ be the anti-linear automorphism of $L_\g(k, 0)$ defined as above. Assume that $\lambda\in\h^*.$ Then $L_\g(k, \lambda) $ is unitary $L_\g(k, 0)$-modules if and only if $\lambda$ is integral dominant and  $(\lambda, \theta)\leq k.$
\end{theorem}
\subsection{Unitary Heisenberg vertex operator algebras}
In this subsection we prove that the Heisenberg vertex operator
algebras are unitary. First, we recall some facts about
Heisenberg vertex operator algebras from \cite{FLM}, \cite{LL}. Let $\h$ be a finite
dimensional vector space of dimension $d$ which has a non-degenerate symmetric bilinear
form $(,)$. Consider the affine algebra
$$\hat{\h}=\h\otimes \C[t, t^{-1}]\oplus \C K$$ with the commutation relations: for $\alpha, \beta\in \h$,
$$[\alpha(m), \beta(n)]=Km(\alpha, \beta)\delta_{m+n,
0},$$$$[\hat{\h}, K]=0,$$
where $\alpha(n)=\alpha\otimes t^n$.

For any $\lambda\in \h$, set $$M_{\h}(1,
\lambda)=U(\hat{\h})/J_{\lambda},$$ where $J_{\lambda}$ is the left
ideal of $U(\hat{\h})$ generated by $\alpha(n),(n\geq 1)$,
$\alpha(0)-(\alpha, \lambda)$ and $K-1$. Set
$e^{\lambda}=1+J_{\lambda}$, we know that $M_{\h}(1, \lambda)$ is
spanned by $\alpha_1(-n_1)\cdots \alpha_k(-n_k)e^{\lambda}$,
$n_1\geq \cdots \geq n_k\geq 1$. Let $\alpha_1,\cdots, \alpha_d$
be an orthonormal basis of $\h$. Set $\omega=\frac{1}{2}\sum_{1\leq
i\leq d}\alpha_i(-1)^21.$ It is well-known that $M_{\h}(1, 0)$
is a vertex operator algebra such that $\1=1+J_0$ is the vacuum vector
and $\omega$ is the Virasoro element \cite{LL}. Furthermore, $M_{\h}(1,
\lambda)$ is an irreducible ordinary $M_{\h}(1, 0)$-module.

In the following we assume that $\h$ is of dimension $1$, i.e.
$\h=\C\alpha$, and that $(\alpha, \alpha)=1$. In this case, we will
denote $M_{\h}(1, 0)$ and $M_{\h}(1, \lambda)$ by $M(1, 0)$ and $M(1,
\lambda)$, respectively. It was proved in \cite{KR} that if
$(\alpha, \lambda) \geq 0$ there exists a unique positive definite Hermitian
form on $M(1, \lambda)$ such that
$$(e^{\lambda}, e^{\lambda})=1,$$$$(\alpha(n)u, v)=(u,
\alpha(-n)v)$$
for $ u, v\in M(1, \lambda).$

\vskip0.25cm
Let $\phi$ be an anti-linear map $\phi:M(1, 0)\to M(1, 0)$
such that
$$\phi(\alpha(-n_1)\cdots \alpha(-n_k))=(-1)^k\alpha(-n_1)\cdots \alpha(-n_k).$$
Note that $\phi \alpha(n)\phi^{-1}=-\alpha(n)$ for $n\in\Z.$ Using a proof similar that of Lemmas \ref{involution}, \ref{involution1} shows that $\phi$ is an anti-linear involution of
$M(1, 0).$

\begin{proposition}
Let $\phi$ be the anti-linear involution of $M(1, 0)$ defined above. Then $(M(1, 0), \phi)$ is a unitary vertex operator algebra and $M(1,
\lambda)$ is a unitary irreducible $M(1, 0)$-module if $(\alpha, \lambda)\geq 0$.
\end{proposition}
\pf By the discussion above, we only need to prove the invariant property. Since $M(1, 0)$ is generated
by $\alpha(-1)$, by Proposition \ref{module} it is
enough to prove that
$$(Y(e^{zL(1)}(-z^{-2})^{L(0)}\alpha(-1), z^{-1})u, v)=(u, Y(\phi(\alpha(-1)), z)v)$$
for $u,v \in M(1, \lambda).$
A straightforward computation gives
 \begin{eqnarray*}
&&(Y(e^{zL(1)}(-z^{-2})^{L(0)}\alpha(-1), z^{-1})u, v)\\
&&=(Y(-z^{-2}\alpha(-1) , z^{-1})u, v)\\
&&=\sum_{n\in \Z}-z^{-2}(\alpha(n)u, v)z^{n+1}\\
&&=\sum_{n\in \Z}(u, -\alpha(-n)v)z^{n-1}\\
&&=(u, Y(\phi(\alpha(-1)), z)v).
\end{eqnarray*}
Then $M(1, \lambda)$ is a unitary module for $M(1,0).$ In particular, $M(1,0)$ is a unitary
vertex operator algebra. \qed

 Note that if $\h$ is a finite dimensional vector space of dimension $d$ and assume that $\alpha_1, ..., \alpha_d$ is an orthonormal basis of $\h$ with respect $(,)$. Then we have $M_{\h}(1, 0)\cong M_{\C\alpha_1}(1,0)\otimes\cdots
\otimes M_{\C\alpha_d}(1,0)$, by Propositions \ref{tensor}, \ref{tensm} we have the following result for general Heisenberg vertex operator algebra $M_{\h}(1, 0)$.
\begin{proposition}\label{heisenberg}
Let $\h$ be a finite dimensional vector space of  dimension $d$ which has a non-degenerate symmetric bilinear
form $(,)$ and $\alpha_1, ..., \alpha_d$ be an orthonormal basis of $\h$ with respect $(,)$. Then $M_{\h}(1, 0)$ is a unitary vertex operator algebra. Furthermore,
if $(\alpha_i, \lambda)\geq 0$, $1\leq i\leq d$, then $M_{\h}(1,
\lambda)$ is a unitary irreducible $M_{\h}(1, 0)$-module.
\end{proposition}
%\begin{remark}
%By Proposition \ref{heisenberg}, we know that the Heisenberg
%vertex operator algebra $M_H(1, 0)$, which is neither rational or
%$C_2$-cofinite, is a unitary vertex operator algebra. But we'll
%see that unitary vertex operator algebras have good property.
%\end{remark}

\subsection{Unitary lattice vertex operator algebras}
In this subsection we prove that the lattice vertex operator
algebras associated to positive definite even lattices are unitary.
First, we recall from \cite{FLM}, \cite{D1} some facts about lattice vertex operator
algebras. Let $L$ be a positive definite even lattice and
$\hat{L}$ be the canonical central extension of $L$ by the cyclic
group $<\kappa>$ of order $2$:
$$1\to <\kappa>\to \hat{L}\to L\to 1,$$ with the commutator map
$c(\alpha, \beta)=\kappa^{(\alpha, \beta)}$. Let $e: L\to \hat{L}$
be a section such that $e_0 = 1$ and $\epsilon_0: L\times L\to <\kappa>$ be the corresponding 2-cocycle. Then $\epsilon_0(\alpha,
\beta)\epsilon_0(\beta, \alpha)=\kappa^{(\alpha, \beta)}$  and
$e_{\alpha}e_{\beta}=\epsilon_0(\alpha, \beta)e_{\alpha+\beta}$ for
$ \alpha, \beta \in L$. Let $\nu:<\kappa>\to <\pm 1>$ be the isomorphism such that $\nu(\kappa)=-1$ and set $$\epsilon=\nu\circ \epsilon_0: L\times L\to <\pm 1>.$$ Consider the induced
$\hat{L}$-module:$$\C\{L\}=\C[\hat{L}]\otimes_{<\kappa>}\C\cong
\C[L]\ \ (linearly),$$ where $\C[.]$ denotes the group algebra and
$\kappa$ acts on $\C$ as multiplication by $-1$.  Then $\C[L]$ becomes a
$\hat{L}$-module such that $e_{\alpha}\cdot
e^{\beta}=\epsilon(\alpha, \beta)e^{\alpha+\beta}$ and $\kappa\cdot
e^{\beta}=-e^{\beta}$. We also define an action $h(0)$ on
$\C[L]$ by: $h(0)\cdot e^{\alpha}=(h, \alpha)e^{\alpha}$ for $h
\in \h, \alpha\in {L}$ and an action $z^h$ on $\C[L]$ by $z^{h}\cdot e^{\alpha}=z^{(h,
\alpha)}e^{\alpha}$.
% For
%$a\in \hat{L},$ write $\iota(a)$ for $a\otimes 1$ in $\C\{L\}$.
%Then the action of $\hat{L}$ on $\C\{L\}$ is given by: $a \cdot
%\iota(b)=\iota(ab)$ and $(-1)\cdot \iota(b)=-\iota(b)$ for
%$\forall a, b\in \hat{L}$.

Set $\h=\C\otimes_{\Z}L$, and  consider the corresponding
Heisenberg vertex operator algebra $M_{\h}(1, 0)$, which is denoted
by $M(1)$ in the following. The untwisted Fock  space associated with
$L$ is defined to be
$$V_L=M(1)\otimes_{\C}\C\{L\}\cong M(1)\otimes_{\C}\C[L]\ (linearly).$$Then $\hat{L}$, $h(n)(n\neq 0)$, $h(0)$ and $z^{h(0)}$ act naturally on $V_L$ by acting on either $M(1)$ or $\C[L]$ as indicated above. It was proved in \cite{B} and \cite{FLM} that
$V_L$ has a vertex operator algebra structure which is determined by
$$Y(h(-1)1, z)=h(z)=\sum_{n\in \Z}h(n)z^{-n-1}\ \ (h\in \h),$$
$$Y(e^{\alpha}, z)=E^-(-\alpha, z)E^+(-\alpha, z)e_{\alpha}z^{\alpha},$$
where
$$E^-(\alpha, z)=exp(\sum_{n<0}\frac{\alpha(n)}{n}z^{-n}),$$
$$E^+(\alpha, z)=exp(\sum_{n>0}\frac{\alpha(n)}{n}z^{-n}).$$

Recall that
$L^{\circ}=\{\,\lambda\in\h\,|\,(\alpha,\lambda)\in\Z\,\}$
is the dual lattice of $L$.  There is an $\widehat{L}$-module
structure on $\C[L^{\circ}]=\bigoplus_{\lambda\in
L^{\circ}}\C e^\lambda$ such that $\kappa$ acts as $-1$
(see \cite{DL}). Let $L^{\circ}=\cup_{i\in
L^{\circ}/L}(L+\lambda_i)$ be the coset decomposition such that
$(\lambda_i,\lambda_i)$ is minimal among all $(\lambda,\lambda)$ for
$\lambda\in L+\lambda_i.$ In particular, $\lambda_0=0.$
Set $\C[L+\lambda_i]=\bigoplus_{\alpha\in L}\C
e^{\alpha+\lambda_i}.$ Then $\C[L^{\circ}]=\bigoplus_{i\in
L^{\circ}/L}\C[L+\lambda_i]$ and each $\C[L+\lambda_i]$ is an
$\widehat L$-submodule of $\C[L^{\circ}].$ The action of
$\widehat L$ on $\C[L+\lambda_i]$ is defined as follow:
$$e_{\alpha}e^{\beta+\lambda_i}=\e(\a,\b)e^{\a+\b+\l_i}$$
for $\alpha,\,\b\in L.$ On the surface, the module structure on each
$\C[L+\lambda_i]$ depends on the choice of $\lambda_i$ in
$L+\lambda_i.$ It is easy to prove that different choices of
$\lambda_i$ give isomorphic $\widehat L$-modules.

 Set $\C[M]=\bigoplus_{\lambda\in M}\C e^{\lambda}$ for a subset $M$ of
$L^{\circ}$, and define $V_M=M(1)\otimes\C[M]$. Then  $V_{L+\lambda_i}$ for $i\in
L^{\circ}/L$ are the irreducible modules for $V_L$ (see \cite{B},
\cite{FLM}, \cite{D1}).

Recall that there is an automorphism $\theta$ of $V_L$ which is
defined as follow: $$\theta:V_L\to V_L$$
$$\alpha_1(-n_1)\cdots\alpha_k(-n_k)\otimes e^{\alpha}\mapsto
(-1)^k\alpha_1(-n_1)\cdots\alpha_k(-n_k)\otimes e^{-\alpha}.$$In particular, $\theta$ induces an automorphism of $M(1)$. Next we recall a construction of $\theta$-twisted modules for $M(1)$ and $V_L$  \cite{FLM}, \cite{D2}.
Denote $\hat{\h}[-1]=\h\otimes t^{1/2}\C[t, t^{-1}]\oplus \C K$ the twisted affinization of $\h$ defined by the communication relations $$[\alpha\otimes t^m, \beta\otimes t^n]=Km(\alpha, \beta)\delta_{m+n,
0},$$$$[\hat{\h}, K]=0,$$for $m, n\in 1/2+\Z$. Then the symmetric algebra $M(1)(\theta)=S(t
^{1/2}\C[t^{-1}]\otimes \h)$ is the unique $\hat{\h}[-1]$-module such that $K=1$ and $\alpha\otimes t^n\cdot 1=0$ if $n>0$. It was proved in \cite{FLM} that $M(1)(\theta)$ is a $\theta$-twisted $M(1)$-module.

By abuse the notation we also use $\theta$ to denote the automorphism of $\hat{L}$ defined by
$\theta(e_{\alpha}) = e_{-\alpha}$ and $\theta(\kappa)=\kappa$. Set
$K=\{\theta(a)a^{-1}|a\in \hat{L}\}$. Let
$\chi$ be a central character of $\hat{L}/K$ such that
$\chi(\kappa K)=-1$ and $T_{\chi}$ be the irreducible
$\hat{L}/K$-module with the character $\chi$.  Then $\hat{L}$, $h(n)$ $(n\in 1/2+\Z)$ act naturally on $V_L^{T_{\chi}}=M(1)(\theta)\otimes T_{\chi}$ by acting on either $M(1)(\theta)$ or $T_{\chi}$ as indicated above. It was proved in \cite{FLM} that
$V_L^{T_{\chi}}=M(1)(\theta)\otimes T_{\chi}$ is an irreducible
$\theta$-twisted $V_L$-module such that
$$Y_{\theta}(\alpha(-1)\cdot 1, z)=\alpha(z)=\sum_{n\in{\frac{1}{2}+\Z}}\alpha(n)z^{-n-1},$$
$$Y_{\theta}(e^{\alpha}, z)=2^{-(\alpha,\alpha)}E^-(-\alpha, z)E^+(-\alpha, z)e_{\alpha}z^{-(\alpha,\alpha)/2},$$
where $$E^{\pm}(\alpha, z)=exp(\sum_{n\in
\pm(\N+\frac{1}{2})}\frac{\alpha(n)}{n}z^{-n}).$$

\vskip0.5cm We now define a Hermitian form on $V_{L^{\circ}}$. First, there is a positive definite Hermitian form on
$\C[L^{\circ}]$:
$$(,):\C[L^{\circ}]\times\C[L^{\circ}]\to \C$$
determined by the conditions  $(e^{\alpha},e^{\beta})=0$ if
$\alpha\neq\beta$ and $(e^{\alpha},e^{\beta})=1$ if
$\alpha=\beta$. And there is
a unique positive definite Hermitian form $(,)$ on $M(1)$ such
that for any $h\in \h$, we have
$$({\bf 1}, {\bf 1})=1,$$$$(h(n)u, v)=(u,
h(-n)v)$$
for all $u, v\in M(1).$

Define a positive definite Hermitian form on $V_{L^{\circ}}$
as follow: for any  $u, v\in M(1)$ and $e^{\alpha}, e^{\beta}\in \C[L^{\circ}]$ $$(u\otimes e^{\alpha}, v\otimes e^{\beta})=(u,
v)(e^{\alpha}, e^{\beta}).$$Note that the positive definite Hermitian form on $V_{L^{\circ}}$ induces a positive definite Hermitian form on $V_{L+\lambda_i}$.

\begin{lem}
Let $(,)$ be the positive definite Hermitian form defined above.
Then we have: for any $\alpha \in L$ and  $ w_1, w_2\in
V_{L^{\circ}}$,
$$(e_{\alpha}w_1, w_2)=(w_1,(-1)^{\frac{(\alpha, \alpha)}{2}} e_{-\alpha}w_2),$$
$$(z^{\alpha}w_1, w_2)=(w_1, z^{\alpha}w_2).$$
\end{lem}
\pf The second identity is obvious. The first identity follows immediately from the fact that $(e_{\alpha} w_1, e_{\alpha}w_2)=(w_1, w_2)$ for any $w_1, w_2\in V_{L^{\circ}}$.\qed

\vskip0.25cm
Let $\phi: V_L\to V_L$ be an anti-linear map which is determined by:
$$\phi: V_L \to V_L$$
$$\alpha_1(-n_1)\cdots\alpha_k(-n_k)\otimes e^{\alpha}\mapsto (-1)^k\alpha_1(-n_1)\cdots\alpha_k(-n_k)\otimes
e^{-\alpha}.$$
Note that we have $$\phi\alpha(n)\phi^{-1}=-\alpha(n)$$
 for $\alpha\in L$ and $n\in \Z,$
and
$$\phi Y(e^{\alpha}, z)\phi^{-1}=Y(e^{-\alpha}, z).$$
Again use a similar discussion as in the proof of Lemmas \ref{involution}, \ref{involution1} shows that $\phi$ is an anti-linear involution of
$V_L.$
\begin{theorem}\label{latt}
Let $L$ be a positive definite even lattice and $\phi$ be the anti-linear involution of $V_L$ defined above. Then the lattice vertex operator algebra $(V_L, \phi)$ is a unitary vertex operator algebra and each
$V_{L+\lambda_i}$ for $i\in
L^{\circ}/L$ is a unitary module for $V_L$.
\end{theorem}
\pf  We only give the proof of the unitarity of $V_L$ here. The proof for $V_{L+\lambda_i}$ is similar.

From the discussion above, we only need to prove the invariant property.  Since  the lattice vertex operator algebra $V_L$ is generated by
$$\{\alpha(-1) | \alpha\in L \}\cup \{e^\alpha |\alpha\in
L\},$$
it is sufficient to prove
the following identities
\begin{equation}\label{1.1}
(Y(e^{zL(1)}(-z^{-2})^{L(0)}\alpha(-1)\cdot1, z^{-1})w_1,
w_2)=(w_1, Y(\phi(\alpha(-1)\cdot 1), z)w_2),\end{equation}
 \begin{equation}\label{1.2}
 (Y(e^{zL(1)}(-z^{-2})^{L(0)}e^{\alpha},
z^{-1})w_1, w_2)=(w_1, Y(\phi(e^{\alpha}), z)w_2)\end{equation}
for any $w_1,w_2 \in V_L$  by Proposition \ref{basic}.

 Assume that
$w_1= u\otimes
e^{\gamma_1},$ $w_2= v\otimes e^{\gamma_2}$ for some $u, v\in M(1)$ and $\gamma_1, \gamma_2\in L $.  By the
definition of the Hermitian form, we have
\begin{eqnarray*}
&&(Y(e^{zL(1)}(-z^{-2})^{L(0)}\alpha(-1), z^{-1})w_1, w_2)\\
&&=-z^{-2}\sum_{n\in \Z}(\alpha(n)w_1, w_2)z^{n+1}\\
&&=\sum_{n\in \Z}-(w_1, \alpha(-n)w_2)z^{n-1}\\
&&=(w_1, Y(\phi(\alpha(-1)), z)w_2).
\end{eqnarray*}
This gives (\ref{1.1}).

To prove the identity (\ref{1.2}), we assume that $(\alpha,
\alpha)=2k$ and $\alpha+\gamma_1= \gamma_2$, then we have
\begin{eqnarray*}
&&(Y(e^{zL(1)}(-z^{-2})^{L(0)}e^{\alpha}, z^{-1})w_1, w_2)\\
&&=(-z^{-2})^k(Y(e^{\alpha}, z^{-1}) u\otimes
e^{\gamma_1}, v\otimes e^{\gamma_2})\\
&&=(-z^{-2})^k(E^-(-\alpha,
z^{-1})E^+(-\alpha, z^{-1})e_{\alpha}(z^{-1})^{\alpha}u\otimes
e^{\gamma_1}, v\otimes
e^{\gamma_2})\\
&&=(-z^{-2})^k(u\otimes
e^{\gamma_1},E^-(\alpha, z)E^+(\alpha,
z)(z^{-1})^{\alpha}(-1)^ke_{-\alpha}v\otimes e^{\gamma_2})\\
&&=(z^{-2})^k(u\otimes
e^{\gamma_1},E^-(\alpha, z)E^+(\alpha,
z)e_{-\alpha}(z^{-1})^{\alpha}(z^{-1})^{(\alpha,-\alpha)}v\otimes e^{\gamma_2})\\
&&=(u\otimes e^{\gamma_1},E^-(\alpha,
z)E^+(\alpha,
z)e_{-\alpha}(z^{-1})^{\alpha}v\otimes e^{\gamma_2})\\
&&=(w_1, Y(e^{-\alpha}, z)w_2)\\
&&=(w_1, Y(\phi(e^{\alpha}), z)w_2).
\end{eqnarray*}
Then $(V_L, \phi)$ is a unitary vertex operator algebra.\qed

We now prove that the $\theta$-twisted $V_L$-module
$V_L^{T_{\chi}}$ is unitary. First, recall from \cite{FLM} that there exists a maximal abelian subgroup  $\hat{\Phi}$ of $\hat{L}$  and
 a homomorphism $\psi:\hat{\Phi}/K\to \C$ extending
$\chi$ such that
$T_{\chi}$ has  form
$Ind_{\hat{\Phi}}^{\hat{L}}\C_{\psi}=\C[\hat{L}]\otimes_{\C[\hat{\Phi}]}\C_{\psi}$. And there is a positive definite Hermitian form $(,):
T_{\chi}\times T_{\chi}\to \C$ on $T_{\chi}$ defined by the
conditions: for any $ a, b \in \hat{L}$, $(t(a), t(b))=0$ if
$a\hat{\Phi}\neq b\hat{\Phi}$  and $(t(a), t(b))=1$ if
$a\hat{\Phi}=b\hat{\Phi}$, where $t(a)=a\otimes 1\in T_{\chi}$ for
$a\in \hat{L}$. Also recall from \cite{FLM} that there is a
positive definite Hermitian form $(,)$ on $M(1)(\theta)$ such that
$$(1,1)=1,$$
$$(h(n)\cdot u, v)=(u, h(-n)\cdot v),$$ for any $u, v\in
M(1)(\theta)$, $h\in \h,$ $n\in 1/2+\Z$. Now we define a positive
definite Hermitian form on $V_L^{T_{\chi}}$ by $(v_1\otimes w_1,
v_2\otimes w_2)=(v_1, v_2)(w_1, w_2)$ for any $v_1, v_2\in
M(1)(\theta)$, $w_1, w_2\in T_{\chi}$.
\begin{lem}\label{1}
Let $(,)$ be the positive definite Hermitian form defined above.
Then we have: for any $\alpha \in L$ and  $ u, v\in
V_L^{T_{\chi}}$,
$$(e_{\alpha}u, v)=(u,(-1)^{\frac{(\alpha, \alpha)}{2}}e_{-\alpha}v). $$
\end{lem}
\pf The lemma follows immediately from the fact that $(e_{\alpha}u, e_{\alpha}v)=(u, v)$ for any $u, v\in V_L^{T_{\chi}}$. \qed
\vskip0.25cm
We now give the unitarity of the $\theta$-twisted $V_L$-module $V_L^{T_{\chi}}$.
\begin{theorem}\label{twist}
For any central character $\chi$, $V_L^{T_{\chi}}$ is a unitary
$\theta$-twisted $V_L$-module.
\end{theorem}
\pf As before, we only need to verify the invariant property in the Definition \ref{dmodule}.  By Proposition \ref{module}, it is sufficient to check
\begin{equation*} (Y(e^{zL(1)}(-z^{-2})^{L(0)}x,
z^{-1})u, v)=(u, Y(\phi(x), z)v)
 \end{equation*}
 for $ x\in \{\alpha(-1)|\alpha \in L\}\cup\{e^{\alpha}|\alpha\in L\}$, $ u, v\in V_L^{T_{\chi}}.$

Assume that $u=v_1\otimes
t(a)$ and $v= v_2\otimes t(b)$ for some $v_1, v_2\in M(1)(\theta)$, $a, b\in \hat{L}$. Then
 $$(\alpha(n) u, v)=(u, \alpha(-n) v)$$ for any $\alpha\in L$ and $n\in 1/2+\Z$. Thus for $x=\alpha(-1)$, we have
\begin{eqnarray*}
&&(Y(e^{zL(1)}(-z^{-2})^{L(0)}\alpha(-1), z^{-1})u, v)\\
&&=-z^{-2}(Y(\alpha(-1), z^{-1})v_1\otimes t(a), v_2\otimes t(b))\\
&&=-z^{-2}\sum_{n\in
\Z+\frac{1}{2}}(\alpha(n)v_1,v_2)(t(a), t(b))z^{n+1}\\
&&=-\sum_{n\in
\Z+\frac{1}{2}}(v_1,\alpha(-n)v_2)(t(a), t(b))z^{n-1}\\
&&=(u, Y(\phi(\alpha(-1)), z)v).
\end{eqnarray*}

Now take $x=e^{\alpha}$ and $(\alpha, \alpha)=2k$. Then by
Lemma \ref{1}  we have
\begin{eqnarray*}
&&(Y(e^{zL(1)}(-z^{-2})^{L(0)}e^{\alpha}, z^{-1})u, v)\\
&&=(Y(e^{zL(1)}(-z^{-2})^{L(0)}e^{\alpha},
z^{-1})v_1\otimes t(a), v_2\otimes t(b))\\
&&=(-z^{-2})^k(2^{-2k}E^-(-\alpha,
z^{-1})E^+(-\alpha, z^{-1})e_{\alpha}z^{k}v_1\otimes t(a),
v_2\otimes t(b))\\
&&=(-z^{-2})^k(v_1\otimes t(a),
2^{-2k}E^-(\alpha,
z)E^+(\alpha, z)(-1)^ke_{-\alpha}z^{k}v_2\otimes t(b))\\
&&=(v_1\otimes t(a), 2^{-2k}E^-(\alpha,
z)E^+(\alpha, z)e_{-\alpha}z^{-k}v_2\otimes t(b))\\
&&=(v_1\otimes t(a), Y(\phi(e^{\alpha}),z)v_2\otimes t(b)).
\end{eqnarray*}
Thus $V_L^{T_{\chi}}$ is a unitary
$\theta$-twisted $V_L$-module.\qed
\subsection{Moonshine vertex operator algebra is unitary}
In this subsection we prove that the famous Moonshine
vertex operator algebra $V^{\natural}$ \cite{FLM} is a unitary vertex
operator algebra. First, we recall some facts about Moonshine
vertex operator algebra. Let $\Lambda$ be the
Leech lattice which is the unique even unimodular lattice with
rank 24 such that $\Lambda_2=\emptyset$. Let $V_{\Lambda}$ be the
lattice vertex operator algebra associated to $\Lambda$ and
$\theta$ be the automorphism of $V_{\Lambda}$ defined as above. Let
$V_{\Lambda}^+=V_{\Lambda}^{\theta}$ be the vertex operator
subalgebra of $V_{\Lambda}$. It was proved in \cite{D3} and \cite{DLM3} that $V_{\Lambda}^+$ is a
rational and $C_2$-cofinite vertex operator algebra. It is
well-known that $V_{\Lambda}^{T}$ is the unique  irreducible $\theta$-twisted
$V_{\Lambda}$-module $V_{\Lambda}^{T}$ \cite{D2} and there is an isomorphism
$\theta$ of $V_{\Lambda}^T$ such that $\theta Y(a,z)\theta^{-1}=Y(\theta(a),z)$ for $a\in V_{\Lambda}$
where $Y(a,z)$ is the twisted vertex operator acting on $V_{\Lambda}^{T}$ \cite{FLM}.
 Let
$(V_{\Lambda}^T)^+=(V_{\Lambda}^T)^{\theta}$, then
$(V_{\Lambda}^T)^+$ is an irreducible $V_{\Lambda}^+$-module. It was proved in \cite{FLM} that $V^{\natural}=V_{\Lambda}^+\oplus
(V_{\Lambda}^T)^+$ is a vertex operator algebra. Furthermore, by the
fusion rule of $V_{\Lambda}^+$ \cite{ADL}, we know that
$(V_{\Lambda}^T)^+$ is a simple current. By Theorems
\ref{latt}, \ref{twist} we  define a positive
definite Hermitian form on $V^{\natural}$ by Theorem
\ref{extension}.
\vskip0.25cm
Define an anti-linear map $\Psi: V_{\Lambda}^T\to
V_{\Lambda}^T$ as follow:
$$\Psi: V_{\Lambda}^T\to V_{\Lambda}^T$$
$$\alpha_1(-n_1)\cdots \alpha_k(-n_k)\otimes t(a)\mapsto (-1)^{k+1}\alpha_1(-n_1)\cdots \alpha_k(-n_k)\otimes t(a).$$
By the definition, $\Psi$ is an anti-linear isomorphism of
$V_{\Lambda}^T$ satisfying
$$\Psi\alpha(n)\Psi^{-1}=-\alpha(n),$$
$$\Psi Y(e^{\alpha}, z)\Psi^{-1}=Y(e^{-\alpha}, z)$$
for $\alpha\in\Lambda.$
That is, $$\Psi Y(\alpha(-1),
z)\Psi^{-1}=Y(\phi(\alpha(-1)), z),$$
$$\Psi Y(e^{\alpha}, z)\Psi^{-1}=Y(\phi(e^{\alpha}), z).$$
This implies
 $$\Psi Y(v, z)\Psi^{-1}=Y(\phi(v), z)$$
  for all $v\in
V_{\Lambda}.$ Note that $\Psi$ commutes with the automorphism  $\theta$ of $V_{\Lambda}^T$ from the definition
of $\theta$ \cite{FLM}. As a result, $\Psi((V_{\Lambda}^T)^+)\subseteq(V_{\Lambda}^T)^+$ and
$\Psi$ induces an anti-linear isomorphism $\psi$
of $(V_{\Lambda}^T)^+$ such that $$\psi Y(v, z)\psi^{-1}=Y(\phi(v),
z)$$
for $v\in V_{\Lambda}^+.$

It is easy to prove that $(\Psi(u), \Psi(v))=\overline{(u,
v)}$ for any $u,v\in V_{\Lambda}^T$. Similarly, for the anti-linear
involution $\phi$ of $V_{\Lambda}$ we could prove that $\phi$ induces an anti-linear involution $\phi_{V_{\Lambda}^+}$ of $V_{\Lambda}^+$ and
$(\phi(x),\phi(y))=\overline{(x, y)}$ for any $x, y \in
V_{\Lambda}$. Then by Theorem \ref{extension} and
Corollary \ref{obi}, we  define an anti-linear involution $\phi_{V^{\natural}}$ of $V^{\natural}$. In fact, we have proved the following:
\begin{theorem}
The Moonshine vertex operator algebra $(V^{\natural}, \phi_{V^{\natural}})$ is a unitary vertex operator algebra.
\end{theorem}
\section{Classification of unitary vertex operator algebras with central charge $c\leq 1$}
\def\theequation{5.\arabic{equation}}
\setcounter{equation}{0}

In this section we consider the classification of unitary
vertex operator algebras with central charge $c\leq 1$. First, we
have the following results about unitary vertex operator algebra with the central charge $c<1$. The similar results about the classification of  rational and $C_2$-cofinite vertex operator algebra  with central charge $c<1$ were obtained in \cite{DZ}, \cite{M}. A classification of local conformal nets with $c<1$ was given in \cite{KL}.
\begin{proposition}
Let $(V, \phi)$ be a unitary vertex operator algebra with central charge
$c<1$. Then the vertex operator subalgebra $<\w>$ which is
generated by $\w$ is isomorphic to the vertex operator algebra $L(c_m, 0)$ for some integer $m\geq 2$. In particular, $V$ is an extension vertex operator algebra of $L(c_m, 0)$.
\end{proposition}
\pf By Lemma \ref{unitary}, we have $V$ is a unitary module of Virasoro
algebra. Since the central charge $c<1$, we have $c=c_m,$
for some integer $m\geq 2.$ By Lemma \ref{unitary}, $V$ is a completely reducible module for
the Virasoro algebra. This implies that the vertex operator
subalgebra of $V$ generated by $\w$ is isomorphic to $L(c_m,0)$ and $V$ is an extension of $L(c_m,0).$
\qed

\vskip0.25cm We now consider the unitary vertex
operator algebra with central charge $c=1$. Recall from \cite{DM2} that for
a vertex operator algebra $V$ of CFT-type, the weight one subspace is
 a Lie algebra under the bracket operation $[a, b]=a_0b$.
\begin{proposition}\label{plast}
Let $(V, \phi)$ be a unitary simple vertex operator algebra of CFT-type with central
charge $c=1$, and  $\h=\C\alpha$ be a one dimensional abelian subalgebra of
$V_1$ such that $\phi(\alpha)\in\R\alpha.$  Then $V$ is
isomorphic to $M_{\h}(1, 0)$ or to $V_L$ for some positive definite even lattice $L$
with rank $1$.
\end{proposition}
\pf We normalize the   Hermitian form $(,)$ so that
$(\1, \1)=1$, then we have $(u, v)=(u_{-1}\1, v)=-(\1,
\phi(u)_1v),$ for any $u, v\in V_1$.

From the assumption we know that $\phi(\alpha)=\lambda\alpha$ for some real number $\lambda.$ Since $\phi$ is an anti-linear involution we see that $\lambda=\pm 1.$ Replacing $\alpha$ by $i\alpha$ is necessary, we can assume that $\phi(\alpha)=-\alpha.$  We can also assume that $(\alpha,\alpha)=1.$ Then
$$[\alpha_m,\alpha_n]=m\delta_{m+n, 0}.$$  By the
Stone-von-Veumann theorem, we have $$V=M_{\h}(1, 0)\otimes
\Omega_V,$$ where $\Omega_V=\{v\in V| \alpha_mv=0,  m>0\}$. Note that the Hermitian form $(,)|_{V_n}$ is
non-degenerate and $(\alpha_0u, v)=(u, \alpha_0v)$ for $ u,v
\in V_n$. This implies that the eigenvalues of $\alpha_0$ on $V_n$ are real. We claim
that $\alpha_0$ acts semisimplely on $V$. Assume  $v\in V_n$ is an generalized eigenvector of $\alpha_0$ with eigenvalue $\lambda$ but not an eigenvector.  Then there exists $n\geq 1$ such that $(\alpha_0-\lambda)^nv\ne 0$ and $(\alpha_0-\lambda)^mv=0$ for $m>n.$ This gives
$$0<((\alpha_0-\lambda)^nv,(\alpha_0-\lambda)^nv)=(v,(\alpha_0-\lambda)^{2n}v)=0,$$
a contradiction. In particular, $\alpha_0$ acts semisimplely on $\Omega_{V}$. For $\lambda\in \h.$

Let
$\Omega_V(\lambda)=\{u\in V|\alpha_0u=(\alpha, \lambda)u, \}.$
Then $$\Omega_V=\oplus_{\lambda\in
\h}\Omega_V(\lambda)$$ and
$$V=\oplus_{\lambda\in \h}M_{\h}(1,
0)\otimes\Omega_V(\lambda).$$Let $L=\{\lambda\in
\h|\Omega_V(\lambda)\neq 0\}$, since $V$ is a simple vertex
operator algebra, by the similar proof as in Theorem 2 of
\cite{DM2} $L$ is an additive subgroup of $\h$. If $L=0$, then we have $V$ is isomorphic to $M_{\h}(1, 0)$. We now assume that $L\neq 0$ and will prove that $L$ is a positive definite lattice.

Set
$$\w'=\frac{1}{2}(\alpha_{-1})^2\1,$$ and $L'(n)=\w'_{n+1}.$ Then
 $$[L'(m), L'(n)]=(m-n)L'(m+n)+\frac{1}{12}(m^3-m)\delta_{m+n, 0}$$ and $$(L'(n)u,v)=(u, L'(-n)v).$$
 In particular, $\w'$ is a Virasoro vector with central charge $1.$
Let $\w''=\w-\w'$, then $\w''$ is a Virasoro vector with central
charge 0. Now we prove that $\w''=0$, i.e. $\w=\w'$. Let
$L''(n)=\w''_{n+1}$, we have
 \begin{eqnarray*}
 &&(\w'', \w'')=(\w-\w',\w-\w')\\
&&\ \ \ \ =((L(-2)-L'(-2))\1, (L(-2)-L'(-2))\1)\\
&&\ \ \ \ =(\1, (L(2)-L'(2))(L(-2)-L'(-2))\1)\\
&&\ \ \ \ =(\1, L''(2)L''(-2)\1)\\
&&\ \ \ \ =0.
 \end{eqnarray*}
 Since the Hermitian form is positive definite, we have $\w''=0$, i.e. $\w=\w'$.
 In particular, $L(0)=L'(0).$ This implies that $\frac{(\lambda, \lambda)}{2}\in \Z^+$ for $\lambda\in L$.
Thus $L$ is a positive definite even lattice. Since the rank of $L$ is $1$,  Corollary 5.4 of
\cite{DM2} then asserts that $V\cong V_L.$ The proof  is complete. \qed

It is not surprised that $M_{\h}(1,0)$ is a possibility in Proposition \ref{plast} as we do not assume that
$V$ is rational. A classification of conformal nets of central charge 1 has been given in \cite{X} under some assumptions.


\begin{thebibliography}{ABCD}
\bibitem[ADL]{ADL} T. Abe, C. Dong and H. Li, Fusion rules for the vertex operator algebras $M(1)^+$ and $V^+
_L$, {\em Comm. Math. Phys.} {\bf 253} (2005), 171-219.
\bibitem[B]{B} R. Borcherds,  Vertex algebras, Kac-Moody algebras, and the
Monster, {\em Proc. Natl. Acad. Sci. USA} {\bf 83} (1986),
3068-3071.
\bibitem[DGM]{DGM} L. Dolan, P. Goddard and P. Montague, Conformal field theories, representations and lattice constructions, {\em Comm. Math. Phys.} {\bf 179}
(1996), 61-120.
\bibitem[D1]{D1} C. Dong, Vertex algebras associated with even lattices, {\em J. Algebra} {\bf 160} (1993), 245¨C265.
\bibitem[D2]{D2} C. Dong, Twisted modules for vertex algebras associated with even lattices, {\em J. Algebra}  {\bf 165} (1994), 91¨C112.
\bibitem[D3]{D3} C. Dong, Representations of the moonshine module vertex operator algebra, {\em Contemp. Math.} {\bf 175} (1994),
27¨C36.
\bibitem[DL]{DL}  C. Dong and J. Lepowsky, Generalized Vertex
Algebras and Relative Vertex Operators, {\em Progress in Math.} Vol.
112, Birkh\"{a}user, Boston 1993.
\bibitem[DLM1]{DLM1} C. Dong, H. Li and G. Mason,
Regularity of rational vertex operator algebras, {\em  Adv. Math.} {\bf 132} (1997), 148-166.
\bibitem[DLM2]{DLM2} C. Dong, H. Li and G. Mason,
Twisted representations of vertex operator algebras, {\em Math. Ann.}
{\bf  310} (1998), 571--600.
\bibitem[DLM3]{DLM3} C. Dong, H. Li and G. Mason,
Modular invariance of the trace functions in the orbifold theory
and generalized moonshine, {\em Comm. Math. Phys.} {\bf 214}
(2000), 1-56.
%\bibitem[DM1]{DM1} C. Dong and G. Mason, On quantum Galois theory, {\em Duke Math. J.} {\bf86} (1997), 305-321.
\bibitem[DM]{DM2} C. Dong and G. Mason, Rational vertex operator
algebra and the effective central charge, {\em Int. Math. Res.
Not.} {\bf 56} (2004), 2989-3008.
\bibitem[DZ]{DZ} C. Dong and W. Zhang, On classfication of rational vertex operator algebras with central charges less than
1, {\em J. Algebra} {\bf 320} (2006), 86-93.
\bibitem[FHL]{FHL}I. Frenkel, Y. Huang and J. Lepowsky, On axiomatic
approaches to vertx operator algebras and modules, {\em Mem. AMS} {\bf  104}, 1993.
\bibitem[FLM]{FLM}I. Frenkel, J. Lepowsky and A. Meurman, Vertex
operator algebras and the Monster, {\em Pure and Appl. Math.} Vol
134, 1988.
\bibitem[FZ]{FZ} I. Frenkel and Y. Zhu,  Vertex operator algebras associated to representations of affine and Virasoro algebras,
{\em Duke Math. J.} {\bf66} (1992), 123-168.
\bibitem[K]{K} V. Kac, Infinite dimensional Lie algebras,  {\em Cambridge Univ. Press, 3er ed}, 1990.
\bibitem[KR]{KR} V. Kac and A. Raina, Bombay lectures on highest
weight representations of infinite dimensional Lie algebras, {\em
World Scientific, Singapore} (1987).
\bibitem[KL]{KL} Y. Kawahigashi and R. Longo, Classification of local conformal nets: Case $c < 1,$ {\em Ann. of Math.} {\bf 160} (2004), 493-522.
\bibitem[LY]{LY} C. Lam and H. Yamauchi, On the structure of framed vertex operator algebras and their pointwise frame
stabilizers, {\em Comm. Math. Phys.} {\bf 277} (2008), 237-285.

\bibitem[LL]{LL} J. Lepowsky and H. Li, Introduction to vertex operator algebras and their representations, {\em Progress in Math.}  Vol 227, 2004.
\bibitem[Li]{Li} H. Li, Symmetric invariant bilinear forms on vertex operator
algebras, {\em Pure and Appl. Math. } {\bf 96} (1994), 279-297.

\bibitem[M]{M} G. Mason,  Lattice subalgebras of strongly regular vertex operator
algebras, {\em arXiv 1110.0544.}
\bibitem[W]{W} W. Wang, Rationality of Virasoro vertex operator algebras,
{\em Internat. Math. Res. Notices} {\bf  7} (1993), 197-211.
\bibitem[X]{X}F. Xu, Strong additivity and conformal nets,
{\em Pacific J. Math.} {\bf 221} (2005), 167-199.
\bibitem[Y]{Y} H. Yamauchi, Module category of simple current extensions of vertex operator algebras, {\em J. Pure Appl. Algebra} {\bf 189} (2004), 315-328.
\bibitem[Z]{Z}Y. Zhu, Modular invariance of characters of vertex
operator algebras, {\em J. AMS} {\bf 9} (1996), 237-302.


\end{thebibliography}
\end{document}